\documentclass[12pt]{article}
\usepackage[utf8]{inputenc}
\usepackage[T1]{fontenc}

\usepackage[]{amsfonts}
\usepackage{amssymb}
\usepackage{amsmath,amsthm}
\def\couleur(#1 #2 #3)
	{
\special{ps::[end]}}

\def\bx#1{\setbox1=\hbox{\kern3pt{#1}\kern3pt}			
 \dimen1=\ht1 \advance\dimen1 by 3pt \dimen2=\dp1 \advance\dimen2 by 3pt
 \setbox1=\hbox{\vrule height\dimen1 depth\dimen2\box1\vrule}%
 \setbox1=\vbox{\hrule\box1\hrule}%
 \advance\dimen1 by .4pt \ht1=\dimen1
 \advance\dimen2 by .4pt \dp1=\dimen2 \box1\relax}

\def\wbb#1{\kern#1em}
\def\vci{\vrule  width.02em height1.47ex depth-.0ex}		
\def\11{{\rm\wbb{.2}\vci\wbb{-.37}1}}

\def\underset#1#2{\mathrel{\mathop{\kern0pt #2}\limits_{#1}}}

\def\overset#1#2{\mathrel{\mathop{\kern0pt #2}\limits^{#1}}}

\newtheorem{thm}{Theorem}[section]
\newtheorem{lem}[thm]{Lemma}
\newtheorem{prop}[thm]{Proposition}
\newtheorem{cor}[thm]{Corollary}
\newtheorem{defin}[thm]{Definition}
\newtheorem{rem}[thm]{Remark}

\begin{document}

\title{Gradient estimates for the heat semigroup on forms in a complete Riemannian  manifold.}

\author{Eric Amar}

\date{ }
\maketitle
 \renewcommand{\abstractname}{Abstract}

\begin{abstract}
We study the heat equation $\frac{\partial u}{\partial t}-\Delta
 u=0,\ u(x,0)=\omega (x),$ where $\Delta :=dd^{*}+d^{*}d$  is
 the Hodge laplacian and $u(\cdot ,t)$ and $\omega $ are $p$-differential
  forms in the complete Riemannian  manifold $(M,g).$ Under weak
 bounded geometrical assumptions we get estimates on its semigroup
 of the form:\par 
acting on $p$-forms with $p\geq 1$ and $k\geq 0$:\par 
\quad \quad \quad $\displaystyle \forall t\geq 1,\ {\left\Vert{\nabla ^{k}e^{-t\Delta
 _{p}}}\right\Vert}_{L^{r}(M)-L^{r}(M)}\leq c(n,r,k).$\par 
Acting on functions, i.e. with $p=0,$ we get a better result:\par 
\quad \quad \quad $\displaystyle \forall k\geq 1,\ \forall t\geq 1,\ {\left\Vert{\nabla
 ^{k}e^{-t\Delta }}\right\Vert}_{L^{r}(M)-L^{r}(M)}\leq c(n,r,k)t^{-1/2}.$\par 

\end{abstract}
\ \par 
\ \par 
\ \par 
\ \par 
\ \par 
\vskip -80pt\ \par 

\tableofcontents

\section{Introduction.}
\quad In physics the heat equation is a partial differential equation
 that describes how the distribution of heat evolves over time
 in a solid medium.\ \par 
\quad For instance in ${\mathbb{R}}^{n}$ the heat equation is the following:
 $\frac{\partial u}{\partial t}-\Delta u=0,\ u(x,0)=\omega (x),$
 where $\Delta $ is the laplacian in ${\mathbb{R}}^{n}.$ In this
 case we have a solution given by the kernel:\ \par 
\ \par 
\quad \quad \quad $\displaystyle \Phi (x,t):={\left\lbrace{
\begin{matrix}
{\frac{1}{(4\pi t)^{n/2}}e^{-\frac{\left\vert{x}\right\vert ^{2}}{4t}}}&{\
 x\in {\mathbb{R}}^{n},\ t>0}\cr 
{0}&{\ x\in {\mathbb{R}}^{n},\ t\leq 0}\cr 
\end{matrix}
}\right.}$\ \par 
\ \par 
and the solution is: $\displaystyle u(x,t):=\int_{{\mathbb{R}}^{n}}{\omega
 (y)\Phi (x-y,t)dy}.$ \ \par 
\quad Clearly working in ${\mathbb{R}}^{n}$ is not enough: take the
 case of a shuttle coming back through the atmosphere. This can
 be modeled by a Riemannian manifold of dimension two, so we
 are naturally leaded to study the heat equation on Riemannian
 manifolds. And now the geometry enters because we need to deal
 with differential forms and the heat equation with them. In
 ${\mathbb{R}}^{n}$ the Laplacian operates diagonally on $p$-forms
 so we have to deal only with functions. This is not the case
 on Riemannian manifolds and the treatment has to be more delicate.\ \par 
\ \par 
\quad In the following $M:=(M,g)$ will be a ${\mathcal{C}}^{\infty
 }$ smooth connected complete Riemannian  manifold without boundary
 unless otherwise stated. We shall just say "Riemannian manifold"
 to mean it.\ \par 
\quad The study of $L^{r}$ estimates for the solutions of the heat
 equation in a Riemannian  manifold started long time ago. A
 basic work was done by R.S. Strichartz~\cite{Strichartz83}.
 In particular he proved that the heat kernel is a contraction
 on the space of functions in $L^{r}(M)$ for $1\leq r\leq \infty .$\ \par 
\quad The study of general parabolic equations in ${\mathbb{R}}^{n}$
 is also well advanced, see for instance ~\cite{HHH06} and the
 references therein. In the case of parabolic equations in Riemannian
 manifold we can see for instance~\cite{Nistor06} and the references
 therein.\ \par 
\quad In~\cite{SobPar18} we also study parabolic equations in vector
 bundles on Riemannian manifold with mixed time-space Lebesgue
 or Sobolev norm. Here we get pointwise in time estimates, and
 we use essentially the same philosophy as in~\cite{SobPar18}
 to pass from local to global by use of the "admissible balls".\ \par 
\ \par 
\quad So our aim in this work is to get estimates on the covariant
 derivatives of any order of solutions of the heat equation $\frac{\partial
 u}{\partial t}-\Delta u=0,\ u(x,0)=\omega (x),$ where $\Delta
 :=dd^{*}+d^{*}d$ is the Hodge laplacian and $u(x,t)$ and $\omega
 (x)$ are $p$-differential  forms in the Riemannian  manifold
 $M.$ We shall denote $L^{r}(U)$ the space of forms, or derivatives
 of them, in the Lebesgue space $L^{r}(U)$ for a measurable set
 $U\subset M,$ with the same notations as for functions\ \par 
\ \par 
\quad We introduce $(m,\epsilon )$-\emph{admissible balls} $B_{m,\epsilon
 }(x)$ in $(M,g)$ as in~\cite{SobPar18}. These balls are the
 ones defined in the work of Hebey and Herzlich~\cite{HebeyHerzlich97}
 but without asking for the harmonicity of the local coordinates.\ \par 

\begin{defin}
Let~\label{mLIR25} $M$ be a Riemannian manifold and $\displaystyle
 x\in M.$ We shall say that the geodesic ball $\displaystyle
 B(x,R)$ is $(0,\epsilon )$-{\bf admissible} if there is a chart
 $\displaystyle (B(x,R),\varphi )$ such that:\par 
\quad 1) $\displaystyle (1-\epsilon )\delta _{ij}\leq g_{ij}\leq (1+\epsilon
 )\delta _{ij}$ in $\displaystyle B(x,R)$ as bilinear forms,\par 
and it will be $(m,\epsilon )$-{\bf admissible} for $m\geq 1,$
 if, moreover:\par 
\quad 2) $\displaystyle \ \sum_{1\leq \left\vert{\beta }\right\vert
 \leq m}{R^{\left\vert{\beta }\right\vert }\sup \ _{i,j=1,...,n,\
 y\in B_{x}(R)}\left\vert{\partial ^{\beta }g_{ij}(y)}\right\vert
 }\leq \epsilon .$\par 
We shall denote ${\mathcal{A}}_{m}(\epsilon )$ the set of $(m,\epsilon
 )$-admissible balls.
\end{defin}

\begin{defin}
~\label{CL26}Let $\displaystyle x\in M,$ we set $R'(x)=\sup \
 \lbrace R>0::B(x,R)\in {\mathcal{A}}(\epsilon )\rbrace .$ We
 shall say that $\displaystyle R_{\epsilon }(x):=\min \ (1,R'(x)/2)$
 is the $\epsilon $-{\bf admissible radius} at $\displaystyle x.$
\end{defin}
\quad We shall follow a natural path to proceed: first we use  known
 result in ${\mathbb{R}}^{n}$ via the Duhamel formula to get
 precise local estimates on $M,$ then we globalise them.\ \par 
\quad Let $x\in M,\ B:=B(x,R)$ be an $\epsilon $-\emph{admissible ball.}\ \par 
\quad $\bullet $ Using Duhamel formula we first get \emph{local estimates
 for any solutions }$u$\emph{ of}\ \par 
\emph{             }$\ \ \ \ \ \ \ \ \ \ \ \ \ \ \ \ \ \ \ \frac{\partial
 u}{\partial t}-\Delta u=0,\ u(x,0)=\omega (x).$\ \par 
\quad $\bullet $ We suppose now that $\omega \in L^{2}(M)\cap L^{r}(M)$
 and, because there is a \emph{global solution} $u(\cdot ,t)\in
 L^{2}(M)$ such that $\frac{\partial u}{\partial t}-\Delta u=0,\
 u(x,0)=\omega (x),$ this global solution verifies also the local
 estimates.\ \par 
\quad $\bullet $ Using Vitali type covering, plus a weight $w(x)$ coming
 from the $\epsilon $-\emph{admissible radius }$R_{\epsilon },$
 we globalise the result.\ \par 
\ \par 
\quad For $p\geq 0$ let $\Lambda ^{p}(M)$ be the set of ${\mathcal{C}}^{\infty
 }$ smooth $p$-forms in $M.$ We know that\ \par 
\quad \quad \quad $\nabla ^{k}:\Lambda ^{p}\rightarrow \Lambda ^{p}\otimes \underbrace{T^{*}M\otimes
 \cdot \cdot \cdot \otimes T^{*}M}_{k\ times},$\ \par 
for the case of general vector bundle with a metric connection
 instead of just the bundle of $p$-forms, see for instance~\cite{Cantor74}
 or ~\cite[Section 2.3, p. 6]{SobPar18} to have the weights added.
 On this tensor product we have a pointwise modulus which allows
 us to define, with a weight $w$:\ \par 
\quad \quad \quad $\forall u\in \Lambda ^{p}(M),\ {\left\Vert{\nabla ^{k}u}\right\Vert}_{L^{r}(M,w)}^{r}:=\int_{M}{\left\vert{\nabla
 ^{k}u}\right\vert ^{r}wdv}.$\ \par 
For instance in the case of a function $u,$ then $\nabla u$ can
 be seen as the $1$-form $du,$ or as the usual gradient vector.
 We compute $\nabla ^{k}u$ locally in Section~\ref{cSG27}, formula~(\ref{cSG21}).\
 \par 
\quad We shall weakened the usual definition of bounded geometry to
 suit our purpose.\ \par 

\begin{defin}
A Riemannian manifold $M$ has $k$-order  {\bf weak bounded geometry} if:\par 
\quad $\bullet $ the injectivity radius $r_{inj}(x)$ at $x\in M$ is
 bounded below by some constant $i>0$ for any $\displaystyle x\in M$\par 
\quad $\bullet $ for $0\leq j\leq k,$ the covariant derivatives $\nabla
 ^{j}Rc$ of the Ricci curvature tensor are bounded in $L^{\infty }(M)$ norm.
\end{defin}
\quad Now we can state our main theorem.\ \par 

\begin{thm}
~\label{cSG26}Let $M$ be a Riemannian manifold. Let $r\in \lbrack
 1,\infty \rbrack .$ For any $\delta >0,$ there is a $\epsilon
 (\delta )>0$ such that for any $\epsilon \leq \epsilon (\delta
 ),$ for any $k\geq 0$ and any $p$-form $\omega \in L^{r}(M)\cap
 L^{2}(M),$ we have, with $u=e^{t\Delta }\omega ,$ the canonical
 solution of the heat equation:\par 
\quad \quad \quad $\displaystyle \forall t\in (\delta ,1),\ {\left\Vert{\nabla
 ^{k}u(\cdot ,t)}\right\Vert}_{L^{r}(M,\ w)}\leq c(n,r)(\frac{\delta
 }{t-\delta ^{3/2}}){\left\Vert{\omega }\right\Vert}_{L^{r}(M)}.$\par 
For $\omega $ any $p$-form with $p\geq 1$ and any $k\geq 0$:\par 
\quad \quad \quad $\displaystyle \forall t\geq 1,\ {\left\Vert{\nabla ^{k}u(\cdot
 ,t)}\right\Vert}_{L^{r}(M,\ w)}\leq c(n,r){\left\Vert{\omega
 }\right\Vert}_{L^{r}(M)},$\par 
and for any function $\omega ,$ we have the better result for $k\geq 1$:\par 
\quad \quad \quad $\displaystyle \forall t\geq 1,\ {\left\Vert{\nabla ^{k}u(\cdot
 ,t)}\right\Vert}_{L^{r}(M,\ w)}\leq c(n,r)t^{-1/2}{\left\Vert{\omega
 }\right\Vert}_{L^{r}(M)}.$\par 
With the weight $w(x):=R_{\epsilon }(x)^{kr+mr}.$\par 
\quad If $k\leq 1$ we have $R_{\epsilon }(x):=R_{m,\epsilon }(x)$ is
 the admissible radius for the $(m,\epsilon )$-admissible balls,
 with $m=1$ if $p=0$ and $m=2$ if $p\geq 1.$  If $k\geq 2,$ then
 $R_{\epsilon }(x)$ is the admissible radius for the $(k,\epsilon
 )$-admissible balls for $p\geq 1$ and for the $(k-1,\epsilon
 )$-admissible balls for $p=0.$
\end{thm}
\ \par 
\quad To get "classical estimates", i.e. estimates without weights,
 we use ~\cite[Corollary, p. 7] {HebeyHerzlich97} and we prove:\ \par 

\begin{thm}
Let $M$ be a Riemannian manifold. Let $r\in \lbrack 1,\infty
 \rbrack $ and $\omega \in L^{r}(M)\cap L^{2}(M).$ For $k=0,\
 1$ suppose that $(M,g)$ has $1$-order weak bounded geometry
 for $p$-forms with $p\geq 1$ and $0$-order weak bounded geometry
 for functions. For $k\geq 2$ suppose that $M$ has $k$-order
 weak bounded geometry for $p$-forms with $p\geq 1$ and $k-1$-order
 weak bounded geometry for functions.\par 
\quad Then the canonical solution $u:=e^{t\Delta }\omega $ of the heat
 equation is such that, for any $k\geq 0$\!\!\!\! , and with
 $\eta =\eta (n,\epsilon ,i,k)$ given by the Corollary~\ref{CF8}
 of Hebey and Herzlich:\par 
\quad \quad \quad $\displaystyle \forall t\in (\delta ,1),\ {\left\Vert{\nabla
 ^{k}u(\cdot ,t)}\right\Vert}_{L^{r}(M)}\leq c(n,r,\eta )(\frac{\delta
 }{t-\delta ^{3/2}}){\left\Vert{\omega }\right\Vert}_{L^{r}(M)}.$\par 
And for $k\geq 0$ and any $p$-form $\omega $:\par 
\quad \quad \quad $\displaystyle \forall t\geq 1,\ {\left\Vert{\nabla ^{k}u(\cdot
 ,t)}\right\Vert}_{L^{r}(M)}\leq c(n,r,\eta ){\left\Vert{\omega
 }\right\Vert}_{L^{r}(M)},$\par 
For functions we get a better estimate for any $k\geq 1$:\par 
\quad \quad \quad $\displaystyle \forall t\geq 1,\ {\left\Vert{\nabla ^{k}u(\cdot
 ,t)}\right\Vert}_{L^{r}(M)}\leq c(n,r,\eta )t^{-1/2}{\left\Vert{\omega
 }\right\Vert}_{L^{r}(M)}.$
\end{thm}
\quad In order to compare with existing result, we deduce from it:\ \par 

\begin{cor}
Let $M$ be a Riemannian manifold. Let $r\in \lbrack 1,\infty
 \rbrack $ and $\omega \in L^{r}(M)\cap L^{2}(M).$ Suppose that
 $M$ has $1$-order weak bounded geometry for $p$-forms with $p\geq
 1$ and $0$-order weak bounded geometry for functions.\par 
\quad Then we get, for the canonical solution $u:=e^{t\Delta }\omega
 $ of the heat equation:\par 
\quad \quad \quad $\displaystyle \forall t\geq 1,\ {\left\Vert{u(\cdot ,t)}\right\Vert}_{L^{r}(M,\
 w)}\leq c(n,r){\left\Vert{\omega }\right\Vert}_{L^{r}(M)}.$\par 
For $t\geq 1$ and acting on $p$-forms with $p\geq 1$:\par 
\quad \quad \quad $\displaystyle \forall t\geq 1,\ {\left\Vert{\nabla u}\right\Vert}_{L^{r}(M)}\leq
 c(n,r){\left\Vert{\omega }\right\Vert}_{L^{r}(M)},$\par 
and acting on functions:\par 
\quad \quad \quad $\displaystyle \forall t\geq 1,\ {\left\Vert{\nabla u}\right\Vert}_{L^{r}(M)}\leq
 c(n,r)t^{-1/2}{\left\Vert{\omega }\right\Vert}_{L^{r}(M)}.$
\end{cor}
\quad Exponential decay of distribution kernels of resolvents was proved
 in the setting of bounded geometry of any order by Kordyukov~\cite{Kordyukov91}.\
 \par 
\ \par 
\quad Recall the volume doubling property for the manifold $M:$ \ \par 
there exists constants $C,D>0,$ such that
\ \par 
\quad \quad \quad $v(x,\lambda r)m<C\lambda ^{D}v(x,r),\ \forall x\in M,\ \forall
 r>0,\ \forall \lambda \geq 1,\ \ \ \ \ \ \ \ \ \ \ \ \ \ \ \
 \ \ \ \ \ \ {\left({D}\right)}$\ \par 
where $v(x,r)=\mu (B(x,r))$ denotes the volume of the ball $\displaystyle
 B(x,r)$ of center $x$ and radius $r.$\ \par 
\quad Let now $p(t,x,y)$ be the heat kernel on functions (the heat
 kernel of the Laplace-
Beltrami operator $\Delta $)\!\!\!\!
 . Recall the Gaussian upper bound for the manifold $M:$ \ \par 
\quad \quad \quad $\displaystyle p(t,x,y)\leq \frac{C}{v(x,{\sqrt{t}})}\mathrm{e}\mathrm{x}\mathrm{p}(-c\frac{\rho
 ^{2}(x,y)}{t})$ $\forall t>0.$$\ \ \ \ \ \ \ \ \ \ \ \ \ \ \
 \ \ \ \ \ \ \ \ \ \ \ \ \ \ \ \ \ \ (G)$  \ \par 
\quad Then J. Magnez and E-M. Ouhabaz~\cite{MagOuh17} proved:\ \par 

\begin{thm}
Suppose that the manifold M has the volume doubling property
 (D),
 the Gaussian upper bound (G) and $R_{k}^{-}\in \hat K$\!\!\!\!
 . Then\par 
\quad (i) the semi group $(e^{-t{\overrightarrow{\Delta _{k}}}})$ acts
 on $L^{p}(\Lambda ^{k}T^{*}M)$ for all $p\in {\left[{1,\infty
 }\right]}$ and\par 
\quad \quad \quad \quad \quad ${\left\Vert{e^{-t{\overrightarrow{\Delta _{k}}}}}\right\Vert}_{p-p}\leq
 C_{p}(t\mathrm{l}\mathrm{o}\mathrm{g}t)^{\left\vert{\frac{1}{2}-\frac{1}{p}}\right\vert
 ^{\frac{D}{2}}},\ t>e.$\par 
\quad (ii) For all $t\geq 1$ and $p\geq 2$\par 
\quad \quad \quad \quad \quad ${\left\Vert{\nabla e^{-t\Delta }}\right\Vert}_{p-p}\leq C_{p}t^{-\frac{1}{p}}.$\par
 
\quad (iii) There exists $C>0$ such that for all $t>0$ and $x,y\in M$\par 
\quad \quad \quad \quad \quad $\left\vert{{\overrightarrow{p_{k}}}(t,x,y)}\right\vert \leq
 C\frac{(1+t+\frac{\rho ^{2}(x,y)}{t})^{\frac{D}{2}}}{v(x,{\sqrt{t}})^{\frac{1}{2}}v(y,{\sqrt{t}})^{\frac{1}{2}}}\mathrm{e}\mathrm{x}\mathrm{p}(-\frac{\rho
 ^{2}(x,y)}{4t}).$\par 
\quad (iv) There exist $C,c>0$ such that for all $t\geq 1$ and $x,y\in M$\par 
\quad \quad \quad \quad \quad $\left\vert{{\overrightarrow{p_{k}}}(t,x,y)}\right\vert \leq
 C\min (1,\frac{t^{\frac{p}{2}}}{v(x,{\sqrt{t}})})\exp (-c\frac{\rho
 ^{2}(x,y)}{t}).$
\end{thm}
With the notation $\displaystyle \left\vert{{\overrightarrow{p_{k}}}(t,x,y)}\right\vert
 $ for the norm from $\Lambda ^{k}T^{*}_{y}M$ to $\displaystyle
 \Lambda ^{k}T^{*}_{y}M$ of the linear map $\displaystyle {\overrightarrow{p_{k}}}(t,x,y)$
 between these two spaces.\ \par 
\quad Comparing to the result of J. Magnez and E-M. Ouhabaz, they get
 Lebesgue estimates on $k$-forms:\ \par 
\quad \quad \quad $\displaystyle {\left\Vert{e^{-t\Delta _{k}}}\right\Vert}_{L^{p}-L^{p}}\leq
 C_{p}(t\log t)^{\left\vert{\frac{1}{2}-\frac{1}{p}}\right\vert
 ^{\frac{D}{2}}},\ t>e,\ p\in \lbrack 1,\infty \rbrack .$\ \par 
And gradient estimates on \emph{functions}:\ \par 
\quad \quad \quad $\displaystyle {\left\Vert{\nabla e^{-t\Delta }}\right\Vert}_{L^{p}-L^{p}}\leq
 C_{p}t^{-\frac{1}{p}},\ t\geq 1,\ p\geq 2.$\ \par 
\quad Here we need that $M$ has $1$ order weak bounded geometry to
 get gradient estimates on $p$-\emph{forms} and we need that
 $M$ has $0$ order weak bounded geometry to get gradient estimates
 on functions. Under these geometric hypotheses, our estimates
 are better. The methods we use are also completely different.\ \par 

\section{Admissible balls.}

\begin{lem}
The $\epsilon $-admissible radius $R_{\epsilon }(x)$ is  continuous.
\end{lem}
\quad Proof.\ \par 
Let $\displaystyle x,y\in M.\ $We set $R'(x)=\sup \ \lbrace R>0::B(x,R)\in
 {\mathcal{A}}(\epsilon )\rbrace .$ Suppose that $R'(x)>d_{g}(x,y),$
 where $d_{g}(x,y)$ is the Riemannian distance between $x$ and
 $y.$ Consider the ball $B(y,\rho )$ of center $y$ and radius
 $\rho :=R'(x)-d_{g}(x,y).$ This ball is contained in $B(x,R'(x))$
 hence, by definition of $R'(x),$ we have that all the points
 in $B(y,\rho )$ verify the conditions 1) and 2) so, by definition
 of $R'(y),$ we have that\ \par 
\quad \quad \quad $\displaystyle R'(y)\geq R'(x)-d_{g}(x,y).$\ \par 
If $\displaystyle R'(x)\leq d_{g}(x,y)$ this is also true because
 $\displaystyle R'(y)>0.$ Exchanging $x$ and $y$ we get that
 $\left\vert{R'(y)-R'(x)}\right\vert \leq d_{g}(x,y).$\ \par 
\quad Hence $R'(x)$ is $1$-lipschitzian so it is continuous. So the
 $\epsilon $-admissible radius $R_{\epsilon }(x)$ is also continuous.
 $\blacksquare $\ \par 

\begin{rem}
~\label{SC16} Because our admissible ball $B(x,R_{\epsilon }(x))$
 is geodesic, we have that the injectivity radius $r_{inj}(x)$
 always verifies $\displaystyle r_{inj}(x)\geq R_{\epsilon }(x).$
\end{rem}

\begin{lem}
~\label{m7}(Slow variation of the admissible radius) Let $M$
 be a Riemannian manifold. With $R(x)=R_{\epsilon }(x),$ the
 $\epsilon $-admissible radius at $x\in M,\ \forall y\in B(x,R(x))$
 we have $R(x)/2\leq R(y)\leq 2R(x).$
\end{lem}
\quad Proof.\ \par 
Let $x,y\in M$ and $d(x,y)$ the Riemannian  distance on $(M,g).$
 Let $y\in B(x,R(x))$ then $d(x,y)\leq R(x)$ and suppose first
 that $R(x)\geq R(y).$ \ \par 
Then, because $R(x)=R'(x)/2,$ we get $y\in B(x,R'(x)/2)$ hence
 we have $B(y,R'(x)/2)\subset B(x,R'(x)).$ But by the definition
 of $R'(x),$ the ball $\displaystyle B(x,R'(x))$ is admissible
 and this implies that the ball $\displaystyle B(y,R'(x)/2)$
 is also admissible for exactly the same constants and the same
 chart; this implies that $R'(y)\geq R'(x)/2$ hence $R(y)\geq
 R(x)/2,$ so $R(x)\geq R(y)\geq R(x)/2.$\ \par 
\quad If $R(x)\leq R(y)$ then\ \par 
\quad \quad \quad $\displaystyle d(x,y)\leq R(x)\Rightarrow d(x,y)\leq R(y)\Rightarrow
 x\in B(y,R'(y)/2)\Rightarrow B(x,R'(y)/2)\subset B(y,R'(y)).$ \ \par 
Hence the same way as above we get $R(y)\geq R(x)\geq R(y)/2\Rightarrow
 R(y)\leq 2R(x).$ So in any case we proved that\ \par 
\quad \quad \quad $\forall y\in B(x,R(x))$ we have $R(x)/2\leq R(y)\leq 2R(x).$
 $\blacksquare $\ \par 

\begin{lem}
~\label{SB30} The $\epsilon $-admissible balls $B(x,R_{\epsilon
 }(x))$ trivialise the bundle $\Lambda ^{p}$ of $p$-forms.
\end{lem}
\quad Proof.\ \par 
Because if $B(x,R)$ is a $\epsilon $-admissible ball, we have
 by Remark~\ref{SC16} that $R\leq r_{inj}(x).$ Then, one can
 choose a local frame field for $\Lambda ^{p}$ on $\displaystyle
 B(x,R)$ by radial parallel translation, as done in~\cite[Section
 13, p. 86-87]{TaylorGD}, see also~\cite[p. 4, eq. (1.3)]{Nistor06}.
 This means that the $\epsilon $-admissible balls also trivialise
 the bundle $\Lambda ^{p}.$ $\blacksquare $\ \par 

\section{Local estimates.}
\quad In order to have the local result, we choose a $(1,\epsilon )$-admissible
 ball $B(x,R)$ and the associated chart $\varphi :B\rightarrow
 {\mathbb{R}}^{n}$ such that $\varphi (x)=0.$ We shall need to
 compare the laplacian $\Delta $ in ${\mathbb{R}}^{n}$ and the
 image $\Delta _{\varphi }$ by $\varphi $ of the laplacian in
 the Riemannian manifold $M.$ For instance for functions we have:\ \par 
\quad \quad \quad $\displaystyle \Delta _{\varphi }f=\frac{1}{{\sqrt{\mathrm{d}\mathrm{e}\mathrm{t}(g_{ij})}}}\partial
 _{i}(g^{ij}{\sqrt{\mathrm{d}\mathrm{e}\mathrm{t}(g_{ij})}}\partial
 _{j}f).$\ \par 
\quad An easy computation gives:\ \par 
\quad \quad \quad $\displaystyle (\Delta _{\varphi }-\Delta )f=(g^{ij}-\delta ^{ij})\partial
 ^{2}_{ij}f+a^{ij}(g)\partial _{i}g^{ij}\partial _{j}f,$\ \par 
where the coefficients $a^{ij}(g)$ are smooth functions of the
 metric $g.$\ \par 
\quad Using the (1) in the definition~\ref{mLIR25} of the admissible
 ball, we get\ \par 
\quad \quad \quad $\displaystyle \left\vert{(g^{ij}-\delta ^{ij})\partial ^{2}_{ij}f}\right\vert
 \leq \epsilon \left\vert{\partial ^{2}f}\right\vert $\ \par 
and using the (2) we get\ \par 
\quad \quad \quad $\displaystyle \left\vert{a^{ij}(g)\partial _{i}g^{ij}\partial
 _{j}f}\right\vert \leq C\epsilon R^{-1}\left\vert{\partial f}\right\vert
 ,$\ \par 
where the constant $C$ depends only on the metric $g.$ So we get\ \par 
\quad \quad \quad $\displaystyle \left\vert{(\Delta _{\varphi }-\Delta )f}\right\vert
 \leq \epsilon \left\vert{\partial ^{2}f}\right\vert +C\epsilon
 R^{-1}\left\vert{\partial f}\right\vert .$\ \par 
\ \par 
\quad To treat the case of $p$-forms, we shall use the Bochner-Weitzenb\"ock
 formula, but in its explicit form, in order to get the dependency
 in the derivatives of the metric tensor.\ \par 
\quad Precisely for a $p$ form $\alpha ,\ p\geq 1,$ the equation (6)
 p. 109 in~\cite{DeRham73} gives in $M$:\ \par 
\quad \quad \quad $\displaystyle (\Delta \alpha )_{k_{1}...,k_{p}}=-\nabla ^{i}\nabla
 _{i}\alpha _{k_{1}...k_{p}}+\sum_{\nu =1}^{p}{(-1)^{\nu }(\nabla
 _{k_{\nu }}\nabla ^{i}-\nabla ^{i}\nabla _{k_{\nu }})\alpha
 _{ik_{1}...\hat k_{\nu }...k_{p}}}.$\ \par 
As is well known, the covariant derivatives are linear in the
 Christoffel symbols, hence in the first derivatives of the metric
 $g.$ Because we apply twice covariant derivation, second order
 derivatives of the metric tensor appear linearly in the sum,
 so this time we need the ball $B$ to be $(2,\epsilon )$-admissible
 and via the chart $\varphi $ we get, the same way as for functions,
 for the image $f$ in ${\mathbb{R}}^{n}$ of the $p$-form $\alpha
 $ in $M$:\ \par 
\quad \quad \quad $\displaystyle \left\vert{(\Delta _{\varphi }-\Delta )f}\right\vert
 \leq \epsilon \left\vert{\partial ^{2}f}\right\vert +C\epsilon
 R^{-2}\left\vert{\partial f}\right\vert .$\ \par 
So we proved:\ \par 

\begin{lem}
~\label{gC16} Let $x$ be a point in a Riemannian manifold $M.$
 Let $\alpha $ be a $p$-form in $\displaystyle L^{r}(B),$ with
 $B:=B(x,R)$ a $(1,\epsilon )$-admissible ball in $M$ if $p=0$
 and a $(2,\epsilon )$-admissible ball in $M$ if $p\geq 1.$ Let
 $\varphi $ be a chart on $B$ and set $\Delta _{\varphi }$ the
 image by $\varphi $ of the laplacian in $(M,g)$ and $f$ the
 image of $\alpha .$ We have that $\Delta _{\varphi }-\Delta
 $ is a second order differential operator of the form:\par 
\quad \quad \quad $\displaystyle (\Delta _{\varphi }-\Delta )f=\sum_{i,j}{a_{ij}\partial
 ^{2}_{ij}f}+\sum_{i}{b_{i}\partial _{i}f}.$\par 
Moreover we get, for $p=0$ in $\varphi (B)$:\par 
\quad \quad \quad $\displaystyle \sum_{i,j}{\left\vert{a_{ij}}\right\vert }\leq
 \epsilon ,\ \sum_{i}{\left\vert{b_{i}}\right\vert }\leq C\epsilon R^{-1},$\par 
and for $p\geq 1$:\par 
\quad \quad \quad $\displaystyle \sum_{i,j}{\left\vert{a_{ij}}\right\vert }\leq
 \epsilon ,\ \sum_{i}{\left\vert{b_{i}}\right\vert }\leq C\epsilon R^{-2}.$
\end{lem}
\quad Now we shall use the Duhamel's formula as in~\cite[Proposition
 3.15]{Rosenberg97}. But, instead to use it to build a parametrix,
 we use it to compare the heat kernel in $\displaystyle {\mathbb{R}}^{n}$
 and the heat kernel in the manifold $M.$\ \par 

\begin{prop}
(Duhamel's formula) Provided $\displaystyle e^{-t(X+Y)}$ exists, we have\par 
\quad \quad \quad $\displaystyle e^{-t(X+Y)}=e^{-tX}-\int_{0}^{t}{e^{-(t-s)(X+Y)}Ye^{-sX}ds}.$
\end{prop}
\quad We apply it to $\displaystyle X+Y:=\Delta _{\varphi },\ X:=\Delta
 $ hence $\displaystyle Y:=\Delta _{\varphi }-\Delta $ where
 $\Delta _{\varphi }$ is the image by $\varphi $ of the laplacian
 on $M,\ \Delta $ is the laplacian on $\displaystyle {\mathbb{R}}^{n}.$\ \par 
\ \par 
\quad Given operators $\displaystyle A(t),\ B(t)$ on our space, we set\ \par 
\quad \quad \quad $\displaystyle A\ast B:=\int_{0}^{t}{A(t-s)B(s)ds}.$\ \par 
\quad By~\cite[formula (3.17)]{Rosenberg97} we get\ \par 
\quad \quad \quad \quad \begin{equation}  e^{-t\Delta _{\varphi }}=e^{-t\Delta }+e^{-t\Delta
 }\ast \sum_{j=1}^{\infty }{(-1)^{j}(Ye^{-t\Delta })^{\ast j}},\
 \label{CF5}\end{equation}\ \par 
because we shall choose $\epsilon $ small enough to make the
 series converging, as we shall see later on.\ \par 
Now on for $\gamma =(\gamma _{1},...,\gamma _{n})\in {\mathbb{N}}^{n}$
 we set $\displaystyle \partial ^{\gamma }f:=\frac{\partial ^{\left\vert{\gamma
 }\right\vert }f}{\partial ^{\gamma _{1}}x_{1}\cdot \cdot \cdot
 \partial ^{\gamma _{n}}x_{n}}$ and $\left\vert{\gamma }\right\vert
 :=\gamma _{1}+\cdot \cdot \cdot +\gamma _{n}.$\ \par 

\begin{prop}
~\label{HF9} Let $r\in \lbrack 1,\infty \rbrack .$ Let $x\in
 M,$ a Riemannian manifold. With $B:=B(x,R)$ a $(m,\epsilon )$-admissible
 ball in $M,$ and any $\delta \in (0,1),$ there is a $\epsilon
 (\delta )>0$ such that for any $\epsilon \leq \epsilon (\delta
 ),$ if $\omega $ is a $p$-form in $\displaystyle L^{r}(B),$
 then the $p$-form $u_{\varphi }:=e^{-t\Delta _{\varphi }}\omega
 _{\varphi }$ verifies, in ${\mathbb{R}}^{n}$\!\!\!\! , with
 $\forall \gamma \in {\mathbb{N}}^{n},$ $l:=\left\vert{\gamma
 }\right\vert /2$\par 
\quad \quad \quad $\displaystyle \forall t\in (\delta ,1),\ {\left\Vert{\partial
 ^{\gamma }u_{\varphi }}\right\Vert}_{L^{r}(B_{\varphi })}\leq
 c(n,r)\frac{\delta }{t-\delta ^{1+l}}R_{\varphi }^{-m}{\left\Vert{\omega
 _{\varphi }}\right\Vert}_{L^{r}(B_{\varphi })}.$\par 
And\par 
\quad \quad \quad $\displaystyle \forall t\geq 1,\ {\left\Vert{\partial ^{\gamma
 }u_{\varphi }}\right\Vert}_{L^{r}(B_{\varphi })}\leq c(n,r,\delta
 )t^{-l}R_{\varphi }^{-m}{\left\Vert{\omega _{\varphi }}\right\Vert}_{L^{r}(B_{\varphi
 })}$\par 
with $m=1$ if $p=0$ and $m=2$ if $p\geq 1.$ And also $B_{\varphi
 }=\varphi (B),\ \omega _{\varphi }=\varphi ^{*}\omega $ etc...
\end{prop}
\quad The long proof of this proposition is postponed to Appendix 1.\ \par 

\subsection{Sobolev comparison estimates.~\label{cSG27}}

\begin{lem}
~\label{CF4}Let $\displaystyle B(x,R)\in {\mathcal{A}}_{m}(\epsilon
 ).$ We have for the Levi-Civita connection on $M$:\par 
\quad \quad \quad $\forall y\in B(x,R),\ \forall k\leq m\in {\mathbb{N}},\ \ \left\vert{\partial
 ^{k-1}\Gamma ^{i}_{lj}(y)}\right\vert \leq C(n,k)\epsilon R^{-k}.$
\end{lem}
\quad Proof.\ \par 
Let $\displaystyle \Gamma ^{k}_{lj}$ be the Christoffel coefficients
 of the Levi-Civita connection  on the tangent bundle $TM.$ We have\ \par 
\quad \quad \quad \begin{equation}  \Gamma ^{i}_{kj}=\frac{1}{2}g^{il}(\frac{\partial
 g_{kl}}{\partial x^{j}}+\frac{\partial g_{lj}}{\partial x^{k}}-\frac{\partial
 g_{jk}}{\partial x^{l}}).\label{SB34}\end{equation}\ \par 
On $\displaystyle B(x,R)\in {\mathcal{A}}_{m}(\epsilon ),$ we
 have $\displaystyle (1-\epsilon )\delta _{ij}\leq g_{ij}\leq
 (1+\epsilon )\delta _{ij}$ as bilinear forms. Hence\ \par 
\quad \quad \quad $\displaystyle \forall y\in B(x,\ R),\ \ \left\vert{\Gamma ^{i}_{kj}(y)}\right\vert
 \leq \frac{3}{2}(1-\epsilon )^{-1}\sum_{\left\vert{\beta }\right\vert
 =1}{\sup \ _{i,j=1,...,n,}\left\vert{\partial ^{\beta }g_{ij}(y)}\right\vert
 }$\ \par 
in a coordinates chart on $\displaystyle B(x,R).$ We also have,
 by definition~\ref{mLIR25} \ \par 
\quad \quad \quad \begin{equation} \sum_{ 1\leq \left\vert{\beta }\right\vert \leq
 m}{R^{\left\vert{\beta }\right\vert }\sup \ _{i,j=1,...,n,\
 y\in B_{x}(R)}\left\vert{\partial ^{\beta }g_{ij}(y)}\right\vert
 }\leq \epsilon .\label{cSG19}\end{equation}\ \par 
Hence\ \par 
\quad \quad \quad $\displaystyle \forall y\in B(x,\ R),\ \ \left\vert{\Gamma ^{i}_{kj}(y)}\right\vert
 \leq \frac{3}{2}(1-\epsilon )^{-1}\epsilon R^{-1}.$\ \par 
Taking the first derivatives on~(\ref{SB34}) gives:\ \par 
\quad \quad \quad $\displaystyle \partial \Gamma ^{i}_{kj}=\frac{1}{2}\partial
 g^{il}(\frac{\partial g_{kl}}{\partial x^{j}}+\frac{\partial
 g_{lj}}{\partial x^{k}}-\frac{\partial g_{jk}}{\partial x^{l}})+\frac{1}{2}g^{il}\partial
 (\frac{\partial g_{kl}}{\partial x^{j}}+\frac{\partial g_{lj}}{\partial
 x^{k}}-\frac{\partial g_{jk}}{\partial x^{l}})$\ \par 
So\ \par 
\quad \quad \quad $\displaystyle \forall y\in B(x,\ R),\ \ \left\vert{\partial
 \Gamma ^{i}_{kj}(y)}\right\vert \leq C(1-\epsilon )^{-1}(\sum_{\left\vert{\beta
 }\right\vert =1}{\sup \ _{i,j=1,...,n,}\left\vert{\partial ^{\beta
 }g_{ij}(y)}\right\vert ^{2}}+\sum_{\left\vert{\beta }\right\vert
 \leq 2}{\sup \ _{i,j=1,...,n,}\left\vert{\partial ^{\beta }g_{ij}(y)}\right\vert
 }).$\ \par 
This gives, using~(\ref{cSG19}):\ \par 
\quad \quad \quad $\displaystyle \forall y\in B(x,\ R),\ \ \left\vert{\partial
 \Gamma ^{i}_{kj}(y)}\right\vert \leq C(1-\epsilon )^{-1}(\epsilon
 ^{2}R^{-2}+\epsilon R^{-2})$\ \par 
and, because $\epsilon <1,$\ \par 
\quad \quad \quad $\displaystyle \forall y\in B(x,\ R),\ \ \left\vert{\partial
 \Gamma ^{i}_{kj}(y)}\right\vert \leq C(1-\epsilon )^{-1}\epsilon
 R^{-2},$\ \par 
the constant $C$ being independent of $x,R$ and $\epsilon .$
 Taking $\epsilon \leq 1/2,$ we get\ \par 
\quad \quad \quad $\displaystyle \forall y\in B(x,\ R),\ \ \left\vert{\partial
 \Gamma ^{i}_{kj}(y)}\right\vert \leq C\epsilon R^{-2},$\ \par 
again the constant $C$ being independent of $x,R$ and $\epsilon .$\ \par 
\quad Derivating $k$ times the formula~(\ref{SB34}), with $k\leq m,$ gives:\ \par 
\quad \quad \quad $\displaystyle \forall y\in B(x,\ R),\ \ \left\vert{\partial
 ^{k-1}\Gamma ^{i}_{kj}(y)}\right\vert \leq C(n,k)\epsilon R^{-k}.$\ \par 
The proof is complete. $\blacksquare $\ \par 

\begin{lem}
~\label{m3}Let $B(x,R)$ be a $(k,\epsilon )$-admissible ball
 in $M$ and $\varphi \ :\ B(x,R)\rightarrow {\mathbb{R}}^{n}$
 be the admissible chart relative to $B(x,R).$ Set $u_{\varphi
 }:=\varphi ^{*}u,$ then, for any $k\in {\mathbb{N}}$:\par 
\quad \quad \quad $\displaystyle {\left\Vert{\nabla ^{k}u}\right\Vert}_{L^{r}(B(x,R))}\leq
 {\left\Vert{\partial ^{k}u_{\varphi }}\right\Vert}_{L^{r}(B_{e}(0,(1+\epsilon
 )R))}+\epsilon \sum_{j=0}^{k-1}{(C_{k}R^{-j-1}){\left\Vert{\partial
 ^{j}u_{\varphi }}\right\Vert}_{L^{r}(B_{e}(0,(1+\epsilon )R))}}.$\par 
and, with $B_{e}(0,t)$ the euclidean ball in ${\mathbb{R}}^{n}$
 centered at $0$ and of radius $t,$\par 
\quad \quad \quad $\displaystyle {\left\Vert{u_{\varphi }}\right\Vert}_{W^{k,r}(B_{e}(0,(1-\epsilon
 )R))}\leq cR^{-k}{\left\Vert{u}\right\Vert}_{W_{p}^{k,r}(B(x,R))}.$\par 
\quad We also have, for $k=0$ and $B(x,R)$ being $(0,\epsilon )$-admissible:\par 
\quad \quad \quad $\displaystyle \forall u\in L_{p}^{r}(B(x,R)),\ {\left\Vert{u}\right\Vert}_{L_{p}^{r}(B(x,R))}\leq
 (1+C\epsilon ){\left\Vert{u_{\varphi }}\right\Vert}_{L^{r}(\varphi
 (B(x,R)))},$\par 
and\par 
\quad \quad \quad $\displaystyle {\left\Vert{u_{\varphi }}\right\Vert}_{L^{r}(B_{e}(0,(1-\epsilon
 )R))}\leq (1+C\epsilon ){\left\Vert{u}\right\Vert}_{L_{p}^{r}(B(x,R))}.$\par 
The constants $c,\ C$ being independent of $B.$\par 
\quad In the case of a function $u$ on $M,$ we have better results.
 Let $B(x,R)$ be a $(k-1,\epsilon )$-admissible ball in $M$ and
 $\varphi \ :\ B(x,R)\rightarrow {\mathbb{R}}^{n}$  be the admissible
 chart relative to $B(x,R).$ Set $u_{\varphi }:=u\circ \varphi
 ^{-1},$ then for $k\geq 1$:\par 
\quad \quad \quad ${\left\Vert{\nabla ^{k}u}\right\Vert}_{L^{r}(B)}\leq {\left\Vert{\partial
 ^{k}u_{\varphi }}\right\Vert}_{L^{r}(B_{e}(0,(1+\epsilon )R))}+$\par 
\quad \quad \quad \quad \quad \quad \quad \quad \quad \quad $\displaystyle +\epsilon C(R^{-1}{\left\Vert{\partial u_{\varphi
 }}\right\Vert}_{L^{r}(B_{e}(0,(1+\epsilon )R))}+\cdot \cdot
 \cdot +R^{-k+1}{\left\Vert{\partial ^{k-1}u_{\varphi }}\right\Vert}_{L^{r}(B_{e}(0,(1+\epsilon
 )R))}).$\par 
and\par 
\quad \quad \quad $\displaystyle {\left\Vert{u_{\varphi }}\right\Vert}_{W^{k,r}(B_{e}(0,(1-\epsilon
 )R))}\leq cR^{1-k}{\left\Vert{u}\right\Vert}_{W^{k,r}(B(x,R))}.$
\end{lem}
\quad Proof.\ \par 
We have to compare the norms of $u,\ \nabla u,\cdot \cdot \cdot
 ,\ \nabla ^{m}u,\ $ with the corresponding ones for $u_{\varphi
 }:=\varphi ^{*}u$ in ${\mathbb{R}}^{n}.$\ \par 
\quad \textbf{Case of a function}\ \par 
Let us start with the case of a function $u$ on $M.$ In this
 case we have:  $\displaystyle (\nabla u)_{j}:=\partial _{j}u\
 $ in local coordinates, so $\displaystyle \left\vert{\nabla
 u(y)}\right\vert =\left\vert{\partial u_{\varphi }(z)}\right\vert .$\ \par 
The components of $\nabla ^{2}u$ are given by $\displaystyle
 (\nabla ^{2}u)_{ij}=\partial _{ij}u-\Gamma ^{k}_{ij}\partial
 _{k}u,$ where the Christoffel symbols $\displaystyle \Gamma
 ^{k}_{ij}$ are those of the Levi-Civita connection. Now we have
 for $B(x,R)$ a $(1,\epsilon )$-admissible ball: $\displaystyle
 \left\vert{\Gamma ^{k}_{ij}}\right\vert \leq C\epsilon /R.$\ \par 
So we get, with $\displaystyle \left\vert{\partial u_{\varphi
 }(z)}\right\vert :=\sum_{j}{\left\vert{\partial _{j}u_{\varphi
 }}\right\vert }$ and $\displaystyle \left\vert{\partial ^{2}u_{\varphi
 }(z)}\right\vert :=\sum_{j,k}{\left\vert{\partial ^{2}_{jk}u_{\varphi
 }}\right\vert },$\ \par 
\quad \quad \quad $\displaystyle \left\vert{\nabla ^{2}u(y)}\right\vert \leq \left\vert{\partial
 ^{2}u_{\varphi }(z)}\right\vert +c\frac{\epsilon }{R}\left\vert{\partial
 u_{\varphi }(z)}\right\vert .$\ \par 
Hence, taking the $L^{r}$ norm, we get\ \par 
\quad \quad \quad $\displaystyle {\left\Vert{\nabla ^{2}u}\right\Vert}_{L^{r}(B(x,R))}\leq
 {\left\Vert{\partial ^{2}u_{\varphi }}\right\Vert}_{L^{r}(B_{e}(0,(1+\epsilon
 )R))}+C\frac{\epsilon }{R}{\left\Vert{\partial u_{\varphi }}\right\Vert}_{L^{r}(B_{e}(0,(1+\epsilon
 )R))}.$\ \par 
For controlling $\nabla ^{k}u$ we need only to have $B(x,R)$
 be $(k-1,\epsilon )$-admissible and we get the same way:\ \par 
\quad \quad \quad $\displaystyle \forall y\in B(x,R),\ \left\vert{\nabla ^{k}u(y)}\right\vert
 \leq \left\vert{\partial ^{k}u_{\varphi }(z)}\right\vert +\epsilon
 (C_{1}R^{-1})\left\vert{\partial u_{\varphi }(z)}\right\vert
 +\cdot \cdot \cdot +C_{k-1}R^{1-k}\left\vert{\partial ^{k-1}u_{\varphi
 }(z)}\right\vert )$\ \par 
So, taking the $L^{r}$ norm, for any $k\in {\mathbb{N}}$ with
 $B:=B(x,R)$ be $(k-1,\epsilon )$-admissible:\ \par 
\quad \quad \quad ${\left\Vert{\nabla ^{k}u}\right\Vert}_{L^{r}(B)}\leq {\left\Vert{\partial
 ^{k}u_{\varphi }}\right\Vert}_{L^{r}(B_{e}(0,(1+\epsilon )R))}+$\ \par 
\quad \quad \quad \quad \quad \quad \quad \quad \quad \quad \quad $\displaystyle +\epsilon C(R^{-1}{\left\Vert{\partial u_{\varphi
 }}\right\Vert}_{L^{r}(B_{e}(0,(1+\epsilon )R))}+\cdot \cdot
 \cdot +R^{-k+1}{\left\Vert{\partial ^{k-1}u_{\varphi }}\right\Vert}_{L^{r}(B_{e}(0,(1+\epsilon
 )R))}).$\ \par 
We shall need also the easy reverse estimate:\ \par 
\quad \quad \quad $\displaystyle {\left\Vert{u_{\varphi }}\right\Vert}_{L^{r}(B_{e}(0,(1-\epsilon
 )R))}\leq {\left\Vert{u}\right\Vert}_{L^{r}(B)}.$\ \par 
\ \par 
\quad \textbf{Case of a }$p$\textbf{-form, }$\displaystyle p\geq 1.$\ \par 
By Lemma~\ref{SB30} the $(m,\epsilon )$-admissible ball $\displaystyle
 B(x,R)$ trivialises the bundle $\Lambda ^{p}$ of $p$-forms on
 $M,$ hence the image of a $p$-form in ${\mathbb{R}}^{n}$ is
 just a vector of functions. Precisely $u_{\varphi }:=\varphi
 ^{*}u\in \varphi (B(x,R)){\times}{\mathbb{R}}^{N}.$\ \par 
We have, because $(1-\epsilon )\delta _{ij}\leq g_{ij}\leq (1+\epsilon
 )\delta _{ij}$ in $B(x,R)$:\ \par 
\quad \quad \quad $\displaystyle B_{e}(0,(1-\epsilon )R)\subset \varphi (B(x,R))\subset
 B_{e}(0,(1+\epsilon )R).$\ \par 
Let $u$ be a $p$-form in $M.$ We have that $\nabla u$ depends
 on the first order derivatives of the metric tensor $g.$ Precisely,
 using formula~\cite[\S  26, p. 106]{DeRham73} set $J:=(i_{1},...,i_{p})\in
 {\mathbb{N}}^{p}$\ \par 
\quad \quad \quad $u_{\varphi }:=\sum_{J\in {\mathbb{N}}^{p}}{\alpha _{J}dx^{J}}=\varphi
 ^{*}u$ in the chart $(B,\varphi )$\!\!\!\! , then we have that
 its covariant derivative $\nabla u$ has for components:\ \par 
\quad \quad \quad \begin{equation}  \nabla _{\partial _{i}}\alpha _{J}=\frac{\partial
 \alpha _{J}}{\partial x^{i}}-\sum_{\nu =1}^{p}{\alpha _{i_{1}...i_{\nu
 -1}ki_{\nu +1}...i_{p}}\Gamma _{i_{\nu }i}^{k}},\label{cSG21}\end{equation}\
 \par 
with summation made with respect to the repeated index $k.$\ \par 
\quad By Lemma~\ref{CF4} we get, with the fact that $B(x,R)$ is $(1,\epsilon
 )$-admissible,\ \par 
\quad \quad \quad $\forall y\in B(x,R),\ \forall k\leq m\in {\mathbb{N}},\ \ \left\vert{\Gamma
 ^{i}_{lj}(y)}\right\vert \leq C(n)\epsilon R^{-1},$\ \par 
with $C$ being independent of $B.$\ \par 
Hence\ \par 
\quad \quad \quad $\displaystyle \forall y\in B(x,R),\ \left\vert{u(y)}\right\vert
 =\left\vert{u_{\varphi }(z)}\right\vert ,\ \ \left\vert{\nabla
 u(y)}\right\vert \leq \left\vert{\partial u}\right\vert +\left\vert{\Phi
 }\right\vert ,$\ \par 
where $\Phi $ is given by formula~(\ref{cSG21}) and depends linearly
 on the coefficients of $u$ and linearly on the first order derivatives
 of the metric tensor $g$ via the Christoffel symbols $\Gamma _{ij}^{k}.$\ \par 
So\ \par 
\quad \quad \quad \begin{equation}  \ \left\vert{\nabla u(y)}\right\vert \leq \left\vert{\partial
 u_{\varphi }(z)}\right\vert +C\epsilon R^{-1}\left\vert{u_{\varphi
 }(z)}\right\vert .\label{SC22}\end{equation}\ \par 
Taking the $L^{r}$ norm of this, we get\ \par 
\quad \quad \quad $\displaystyle \ {\left\Vert{\nabla u(y)}\right\Vert}_{L^{r}(B(x,R))}\leq
 {\left\Vert{\partial u_{\varphi }}\right\Vert}_{L^{r}(B_{e}(0,(1+\epsilon
 )R))}+C\epsilon R^{-1}{\left\Vert{u_{\varphi }}\right\Vert}_{L^{r}(B_{e}(0,(1+\epsilon
 )R))}.$\ \par 
The same way for $\nabla ^{k}u$ with $1<k\leq m,$ by iterating
 formula~(\ref{cSG21}) and still with Lemma~\ref{CF4}, we have:\ \par 
\quad \quad \quad $\displaystyle \forall y\in B(x,R),\ \left\vert{\nabla ^{k}u(y)}\right\vert
 \leq \left\vert{\partial ^{k}v(z)}\right\vert +$\ \par 
\quad \quad \quad \quad \quad \quad \quad $\displaystyle +\epsilon (C_{0}R^{-1}\left\vert{u_{\varphi }(z)}\right\vert
 +\ C_{1}R^{-2}\left\vert{\partial u_{\varphi }(z)}\right\vert
 +\cdot \cdot \cdot +C_{k-1}R^{-k}\left\vert{\partial ^{k-1}u_{\varphi
 }(z)}\right\vert ).$\ \par 
We deduce\ \par 
\quad \quad \quad $\displaystyle {\left\Vert{\nabla ^{k}u}\right\Vert}_{L^{r}(B(x,R))}\leq
 {\left\Vert{\partial ^{k}u_{\varphi }}\right\Vert}_{L^{r}(B_{e}(0,(1+\epsilon
 )R))}+\epsilon \sum_{j=0}^{k-1}{(C_{k}R^{-j-1}){\left\Vert{\partial
 ^{j}u_{\varphi }}\right\Vert}_{L^{r}(B_{e}(0,(1+\epsilon )R))}}.$\ \par 
So, with a new constant $c$ independent of $B$:\ \par 
\begin{equation} {\left\Vert{ \nabla ^{k}u}\right\Vert}_{L^{r}(B(x,R))}\leq
 {\left\Vert{\partial ^{k}u_{\varphi }}\right\Vert}_{L^{r}(B_{e}(0,(1+\epsilon
 )R))}+c\epsilon R^{-k}\sum_{j=0}^{k-1}{{\left\Vert{\partial
 ^{j}u_{\varphi }}\right\Vert}_{L^{r}(B_{e}(0,(1+\epsilon )R))}},\label{m2}\end{equation}\
 \par 
because $R\leq 1.$\ \par 
And, using the fact that $\displaystyle {\left\Vert{u}\right\Vert}_{W^{m,r}(B(x,R))}\simeq
 \sum_{k=0}^{m}{{\left\Vert{\nabla ^{k}u}\right\Vert}_{L^{r}(B(x,R))}},$
 we also get:\ \par 
\quad \quad \quad $\displaystyle {\left\Vert{u}\right\Vert}_{W^{m,r}(B(x,R))}\leq
 cR^{-m}{\left\Vert{u_{\varphi }}\right\Vert}_{W^{m,r}(B_{e}(0,(1+\epsilon
 )R))}.$\ \par 
The same way we get the reverse estimates\ \par 
\quad \quad \quad $\displaystyle {\left\Vert{u_{\varphi }}\right\Vert}_{W^{m,r}(B_{e}(0,(1-\epsilon
 )R))}\leq cR^{-m}{\left\Vert{u}\right\Vert}_{W^{m,r}(B(x,R)}.$\ \par 
The case $m=0$ is given by the equation~(\ref{m2}).\ \par 
All the constants here are independent of $B.$\ \par 
\quad The proof of the lemma is complete. $\blacksquare $\ \par 

\subsection{The main local estimates.}

\begin{thm}
~\label{m6}Let $r\in \lbrack 1,\infty \rbrack $ and $k\in {\mathbb{N}}.$
 Let $\displaystyle B:=B(x,R)$ be a $(\beta ,\epsilon )$-admissible
 ball in the Riemannian manifold $M$ with $\beta :=\max (k-1,1)$
 if $p=0$ and $\beta :=\max (k,2)$ if $p\geq 1.$ Then, with $\omega
 \in L_{p}^{2}(B)\cap L_{p}^{r}(B)$ and $u(x,t)=e^{-t\Delta }\omega
 ,$ the canonical solution of the heat equation, we get:\par 
\quad \quad \quad $\displaystyle \forall t\in (\delta ,1),\ {\left\Vert{\nabla
 ^{k}u(\cdot ,t)}\right\Vert}_{L^{r}(B(x,R))}\leq c(n,r)\frac{\delta
 }{t-\delta ^{3/2}}R^{-k-m}{\left\Vert{\omega }\right\Vert}_{L^{r}(B)}.$\par 
and\par 
\quad \quad \quad $\displaystyle \forall t\geq 1,\ {\left\Vert{u(\cdot ,t)}\right\Vert}_{L^{r}(B)}\leq
 c(n,r)R^{-m}{\left\Vert{\omega }\right\Vert}_{L^{r}(B)}.$\par 
This gives for $p$-forms with $p\geq 1,\ k\geq 1,$\par 
\quad \quad \quad $\displaystyle \forall t\geq 1,\ {\left\Vert{\nabla ^{k}u(\cdot
 ,t)}\right\Vert}_{L^{r}(B)}\leq c(n,r)R_{\varphi }^{-k-m}{\left\Vert{\omega
 }\right\Vert}_{L^{r}(B)}.$\par 
And for functions, i.e. $p=0,\ k\geq 1,$\par 
\quad \quad \quad $\displaystyle \forall t\geq 1,\ {\left\Vert{\nabla ^{k}u(\cdot
 ,t)}\right\Vert}_{L^{r}(B)}\leq c(n,r)t^{-1/2}R_{\varphi }^{-k-m}{\left\Vert{\omega
 }\right\Vert}_{L^{r}(B)}.$\par 
Again we have that $m=1$ on functions and $m=2$ on $p$-forms with $p\geq 1.$
\end{thm}
\quad Proof.\ \par 
Let $\omega \in L^{2}(B)\cap L^{r}(B)$ and $u=e^{-t\Delta }\omega
 $ the canonical solution of the heat equation.\ \par 
The ball $B$ being admissible, there is a diffeomorphism $\varphi
 :B\rightarrow {\mathbb{R}}^{n}$ such that $\Lambda ^{p}$ trivialises
 on $B.$\ \par 
So the local representation of the $p$-form $u$ is a vector of functions.\ \par 
We shall apply Proposition~\ref{HF9}, with $m=1$ if $p=0$ and
 $m=2$ if $p\geq 1,$ and with $\displaystyle \forall \gamma \in
 {\mathbb{N}}^{n},\ l:=\left\vert{\gamma }\right\vert /2:$\ \par 
\quad \quad \quad $\displaystyle \forall t\in (\delta ,1),\ \ {\left\Vert{\partial
 ^{\gamma }u_{\varphi }(\cdot ,t)}\right\Vert}_{L^{r}(B_{\varphi
 })}\leq c(n,r)\frac{\delta }{t-\delta ^{1+l}}R_{\varphi }^{-m}{\left\Vert{\omega
 _{\varphi }}\right\Vert}_{L^{r}(B_{\varphi })}.$\ \par 
\ \par 
\quad \quad \quad \begin{equation}  \forall t\geq 1,\ \forall \gamma \in {\mathbb{N}}^{n},\
 {\left\Vert{\partial ^{\gamma }u_{\varphi }(\cdot ,t)}\right\Vert}_{L^{r}(B_{\varphi
 })}\leq c(n,r,\delta )t^{-l}R_{\varphi }^{-m}{\left\Vert{\omega
 _{\varphi }}\right\Vert}_{L^{r}(B_{\varphi })}.\label{cSG23}\end{equation}\
 \par 
Where $B_{\varphi },u_{\varphi },\omega _{\varphi }$ are the
 images by $\varphi $ of $B,u,\omega $ and the image of $\Lambda
 ^{p}$ is the trivial bundle  $\varphi (B){\times}{\mathbb{R}}^{N}$
 in $\displaystyle {\mathbb{R}}^{n}.$ The constants being independent
 of $B.$\ \par 
\ \par 
\quad First, because of the condition $\displaystyle (1-\epsilon )\delta
 _{ij}\leq g_{ij}\leq (1+\epsilon )\delta _{ij}$ in the definition
 of the $\epsilon $-admissible ball, we have that $R_{\varphi
 }\simeq R.$ Recall that $u=e^{-t\Delta }\omega .$\ \par 
\quad Now we use the Sobolev comparison estimates given by Lemma~\ref{m3}
 and, to apply it, we need to have $B\in {\mathcal{A}}_{k}(\epsilon
 )$ and this is the reason to define $\beta :=\max (k-1,1)$ if
 $p=0$ and $\beta :=\max (k,2)$ if $p\geq 1.$ And we get:\ \par 
\quad \quad \quad \begin{equation} {\left\Vert{ u(\cdot ,t)}\right\Vert}_{L^{r}(B)}\leq
 C{\left\Vert{u_{\varphi }(\cdot ,t)}\right\Vert}_{L^{r}(B_{\varphi
 })},\label{cSG20}\end{equation}\ \par 
and, for any $k\in {\mathbb{N}},$\ \par 
\quad \quad \quad $\displaystyle {\left\Vert{\nabla ^{k}u(\cdot ,t)}\right\Vert}_{L^{r}(B(x,R))}\leq
 {\left\Vert{\partial ^{k}u_{\varphi }(\cdot ,t)}\right\Vert}_{L^{r}(B_{\varphi
 })}+\epsilon \sum_{j=0}^{k-1}{(C_{j}R^{-j-1}){\left\Vert{\partial
 ^{j}u_{\varphi }(\cdot ,t)}\right\Vert}_{L^{r}(B_{\varphi })}},$\ \par 
hence, because $R\leq 1,$\ \par 
\quad \quad \quad \begin{equation} {\left\Vert{ \nabla ^{k}u(\cdot ,t)}\right\Vert}_{L^{r}(B(x,R))}\leq
 {\left\Vert{\partial ^{k}u_{\varphi }(\cdot ,t)}\right\Vert}_{L^{r}(B_{\varphi
 })}+\epsilon CR^{-k}{\left\Vert{u_{\varphi }(\cdot ,t)}\right\Vert}_{W^{k-1,r}(B_{\varphi
 })}.\label{nC31}\end{equation}\ \par 
\quad \textbf{Case of functions}\ \par 
In the case of functions, we need only to have $B\in {\mathcal{A}}_{\beta
 }(\epsilon )$ with $\beta :=\max (k-1,1)$ and we have no term
 in $u_{\varphi }$ in the right hand side, so we get:\ \par 
\quad \quad \quad $\displaystyle {\left\Vert{\nabla ^{k}u(\cdot ,t)}\right\Vert}_{L^{r}(B(x,R))}\leq
 {\left\Vert{\partial ^{k}u_{\varphi }(\cdot ,t)}\right\Vert}_{L^{r}(B_{\varphi
 })}+\epsilon \sum_{j=1}^{k-1}{(C_{j}R^{-j-1}){\left\Vert{\partial
 ^{j}u_{\varphi }(\cdot ,t)}\right\Vert}_{L^{r}(B_{\varphi })}},$\ \par 
hence:\ \par 
\quad \quad \quad \begin{equation} {\left\Vert{ \nabla ^{k}u(\cdot ,t)}\right\Vert}_{L^{r}(B(x,R))}\leq
 {\left\Vert{\partial ^{k}u_{\varphi }(\cdot ,t)}\right\Vert}_{L^{r}(B_{\varphi
 })}+\epsilon CR^{-k}\sum_{j=1}^{k-1}{{\left\Vert{\partial ^{j}u_{\varphi
 }(\cdot ,t)}\right\Vert}_{L^{r}(B_{\varphi })}},\label{cSG22}\end{equation}\
 \par 
The constants being independent of $B.$\ \par 
\quad Now, still by Lemma~\ref{m3},\ \par 
\quad \quad \quad $\displaystyle {\left\Vert{\omega _{\varphi }}\right\Vert}_{L^{r}(B_{\varphi
 })^{N}}\leq C{\left\Vert{\omega }\right\Vert}_{L^{r}(B)}.$\ \par 
Hence replacing in~(\ref{cSG20}) we get, with new constants (with
 $m=1$  because we deal with functions) :\ \par 
\quad \quad \quad $\displaystyle \forall t\in (\delta ,1),\ {\left\Vert{u(\cdot
 ,t)}\right\Vert}_{L^{r}(B)}\leq c(n,r){\left\Vert{u_{\varphi
 }(\cdot ,t)}\right\Vert}_{L^{r}(B_{\varphi })}\leq $\ \par 
\quad \quad \quad \quad \quad \quad \quad \quad \quad $\displaystyle \leq c(n,r)\frac{\delta }{t-\delta }R^{-m}{\left\Vert{\omega
 _{\varphi }}\right\Vert}_{L^{r}(B_{\varphi })}\leq c(n,r)\frac{\delta
 }{t-\delta }R^{-m}{\left\Vert{\omega }\right\Vert}_{L^{r}(B)}.$\ \par 
So\ \par 
\quad \quad \quad $\displaystyle \forall t\in (\delta ,1),\ {\left\Vert{u(\cdot
 ,t)}\right\Vert}_{L^{r}(B)}\leq c(n,r)\frac{\delta }{t-\delta
 }R^{-m}{\left\Vert{\omega }\right\Vert}_{L^{r}(B)}.$\ \par 
The same way:\ \par 
\quad \quad \quad $\displaystyle \forall t\geq 1,\ {\left\Vert{u(\cdot ,t)}\right\Vert}_{L^{r}(B)}\leq
 c(n,r)R^{-m}{\left\Vert{\omega }\right\Vert}_{L^{r}(B)}.$\ \par 
\quad For the gradient estimate on functions, we get as above, with
 $\nabla $ the covariant derivative on $M$ and $t\in (\delta ,1):$\ \par 
\quad \quad \quad $\displaystyle {\left\Vert{\nabla u(\cdot ,t)}\right\Vert}_{L^{r}(B(x,R))}\leq
 {\left\Vert{\partial u_{\varphi }(\cdot ,t)}\right\Vert}_{L^{r}(B_{\varphi
 })}\leq c(n,r)\delta ^{-1/2}\frac{\delta ^{3/2}}{t-\delta ^{3/2}}R_{\varphi
 }^{-m}{\left\Vert{\omega _{\varphi }}\right\Vert}_{L^{r}(B_{\varphi
 })}.$\ \par 
Hence\ \par 
\quad \quad \quad $\displaystyle \forall t\in (\delta ,1),\ {\left\Vert{\nabla
 u(\cdot ,t)}\right\Vert}_{L^{r}(B)}\leq c(n,r)\frac{\delta }{t-\delta
 ^{3/2}}R^{-m}{\left\Vert{\omega }\right\Vert}_{L^{r}(B)}.$\ \par 
\ \par 
For functions we have no term in $u_{\varphi }$ by Lemma~\ref{CF4}:\ \par 
\quad \quad \quad $\displaystyle {\left\Vert{\nabla u(\cdot ,t)}\right\Vert}_{L^{r}(B(x,R))}\leq
 {\left\Vert{\partial u_{\varphi }(\cdot ,t)}\right\Vert}_{L^{r}(B_{\varphi
 })},$\ \par 
hence by~\ref{cSG23}\ \par 
\quad \quad \quad $\displaystyle \forall t\geq 1,\ {\left\Vert{\nabla u(\cdot ,t)}\right\Vert}_{L^{r}(B)}\leq
 c(n,r)t^{-1/2}R^{-m}{\left\Vert{\omega }\right\Vert}_{L^{r}(B)}.$\ \par 
\quad And more generally, for $k\geq 2,$ by the same way, with this
 time $B(x,R)$ being a $(k-1,\epsilon )$-admissible ball for functions:\ \par 
\quad \quad \quad $\displaystyle \forall t\in (\delta ,1),\ {\left\Vert{\nabla
 ^{k}u(\cdot ,t)}\right\Vert}_{L^{r}(B(x,R))}\leq c(n,r)\frac{\delta
 }{t-\delta ^{1+k/2}}R_{\varphi }^{-m}{\left\Vert{\omega _{\varphi
 }}\right\Vert}_{L^{r}(B_{\varphi })}+$\ \par 
\quad \quad \quad \quad \quad \quad \quad \quad \quad \quad \quad \quad \quad \quad \quad \quad \quad \quad \quad \quad \quad $\displaystyle +\epsilon c(n,r)CR^{-k}\sum_{j=1}^{k-1}{\frac{\delta
 }{t-\delta ^{1+j/2}}R_{\varphi }^{-m}{\left\Vert{\omega _{\varphi
 }}\right\Vert}_{L^{r}(B_{\varphi })},}$\ \par 
with $l=\left\vert{\gamma }\right\vert /2$ and where we used\ \par 
\quad \quad \quad $\displaystyle \forall t\in (\delta ,1),\ \ {\left\Vert{\partial
 ^{\gamma }u_{\varphi }(\cdot ,t)}\right\Vert}_{L^{r}(B_{\varphi
 })}\leq c(n,r)\frac{\delta }{t-\delta ^{1+l}}R_{\varphi }^{-m}{\left\Vert{\omega
 _{\varphi }}\right\Vert}_{L^{r}(B_{\varphi })}.$\ \par 
Hence\ \par 
\quad \quad \quad $\displaystyle \forall t\in (\delta ,1),\ {\left\Vert{\nabla
 ^{k}u(\cdot ,t)}\right\Vert}_{L^{r}(B(x,R))}\leq c(n,r)\frac{\delta
 }{t-\delta ^{3/2}}R_{}^{-k-m}{\left\Vert{\omega }\right\Vert}_{L^{r}(B)}.$\
 \par 
\quad For $t\geq 1$ we have by ~\ref{cSG22}:\ \par 
\quad \quad \quad $\displaystyle \forall t\geq 1,\ {\left\Vert{\nabla ^{k}u(\cdot
 ,t)}\right\Vert}_{L^{r}(B(x,R))}\leq {\left\Vert{\partial ^{k}u_{\varphi
 }(\cdot ,t)}\right\Vert}_{L^{r}(B_{\varphi })}+\epsilon CR^{-k}\sum_{j=1}^{k-1}{{\left\Vert{\partial
 ^{j}u_{\varphi }(\cdot ,t)}\right\Vert}_{L^{r}(B_{\varphi })}},$\ \par 
so\ \par 
\quad \quad \quad $\displaystyle \forall t\geq 1,\ {\left\Vert{\nabla ^{k}u(\cdot
 ,t)}\right\Vert}_{L^{r}(B)}\leq c(n,r)t^{-1/2}R_{\varphi }^{-k-m}{\left\Vert{\omega
 }\right\Vert}_{L^{r}(B)}.$\ \par 
\ \par 
\quad \textbf{Case of forms}\ \par 
This time $B(x,R)$ is a $(k,\epsilon )$-admissible ball (with
 $m=2$  because we deal with $p$-forms with $p\geq 1$).\ \par 
\quad For forms because by~(\ref{cSG23}) \ \par 
\quad \quad \quad $\displaystyle \forall t\geq 1,\ \forall \gamma \in {\mathbb{N}}^{n},\
 l:=\left\vert{\gamma }\right\vert /2:\ \ {\left\Vert{\partial
 ^{\gamma }u_{\varphi }(\cdot ,t)}\right\Vert}_{L^{r}(B_{\varphi
 })}\leq c(n,r,\delta )t^{-l}R_{\varphi }^{-m}{\left\Vert{\omega
 _{\varphi }}\right\Vert}_{L^{r}(B_{\varphi })}.$\ \par 
Hence\ \par 
\quad \quad \quad $\displaystyle \forall t\geq 1,\ {\left\Vert{\nabla u_{\varphi
 }(\cdot ,t)}\right\Vert}_{L^{r}(B)}\leq c(n,r)t^{-1/2}R_{\varphi
 }^{-m}{\left\Vert{\omega _{\varphi }}\right\Vert}_{L^{r}(B_{\varphi
 })}$ \ \par 
and putting it in~\ref{nC31} we get:\ \par 
\quad \quad \quad $\displaystyle \forall t\geq 1,\ {\left\Vert{\nabla u(\cdot ,t)}\right\Vert}_{L^{r}(B)}\leq
 c(n,r)(t^{-1/2}R_{\varphi }^{-m}{\left\Vert{\omega _{\varphi
 }}\right\Vert}_{L^{r}(B_{\varphi })}+\epsilon CR_{\varphi }^{-1}{\left\Vert{u_{\varphi
 }}\right\Vert}_{L^{r}(B_{\varphi })}$\ \par 
hence with $\gamma =0$ in ~\ref{cSG23}\ \par 
\quad \quad \quad $\displaystyle \forall t\geq 1,\ \ {\left\Vert{u_{\varphi }(\cdot
 ,t)}\right\Vert}_{L^{r}(B_{\varphi })}\leq c(n,r,\delta )R_{\varphi
 }^{-m}{\left\Vert{\omega _{\varphi }}\right\Vert}_{L^{r}(B_{\varphi
 })},$\ \par 
so\ \par 
\quad \quad \quad $\displaystyle \forall t\geq 1,\ {\left\Vert{\nabla u(\cdot ,t)}\right\Vert}_{L^{r}(B)}\leq
 c(n,r)t^{-1/2}R_{\varphi }^{-m}{\left\Vert{\omega _{\varphi
 }}\right\Vert}_{L^{r}(B_{\varphi })}+$\ \par 
\quad \quad \quad \quad \quad \quad \quad \quad \quad \quad \quad \quad \quad \quad \quad $\displaystyle +\epsilon CR_{\varphi }^{-1}c(n,r,\delta )R_{\varphi
 }^{-m}{\left\Vert{\omega _{\varphi }}\right\Vert}_{L^{r}(B_{\varphi
 })},$\ \par 
hence\ \par 
\quad \quad \quad $\displaystyle \forall t\geq 1,\ {\left\Vert{\nabla u(\cdot ,t)}\right\Vert}_{L^{r}(B)}\leq
 c(n,r,\delta )R_{\varphi }^{-m}{\left\Vert{\omega _{\varphi
 }}\right\Vert}_{L^{r}(B_{\varphi })}(t^{-1/2}+\epsilon CR_{\varphi
 }^{-1}).$\ \par 
And by use of Lemma~\ref{m3} and because for any $t\geq 1,\ t^{-1/2}\leq
 1$ and $R<1,$\ \par 
\quad \quad \quad $\displaystyle \forall t\geq 1,\ {\left\Vert{\nabla u(\cdot ,t)}\right\Vert}_{L^{r}(B)}\leq
 c(n,r,\delta )R^{-1-m}{\left\Vert{\omega }\right\Vert}_{L^{r}(B)}.$\ \par 
So for $p$-forms with $p\geq 1,$\ \par 
\quad \quad \quad \begin{equation} {\left\Vert{ \nabla ^{k}u(\cdot ,t)}\right\Vert}_{L^{r}(B(x,R))}\leq
 {\left\Vert{\partial ^{k}u_{\varphi }(\cdot ,t)}\right\Vert}_{L^{r}(B_{\varphi
 })}+\epsilon CR^{-k}{\left\Vert{u_{\varphi }(\cdot ,t)}\right\Vert}_{W^{k-1,r}(B_{\varphi
 })}.\label{cSG24}\end{equation}\ \par 
Hence\ \par 
\quad \quad \quad $\displaystyle \forall t\geq 1,\ {\left\Vert{\nabla ^{k}u(\cdot
 ,t)}\right\Vert}_{L^{r}(B)}\leq c(n,r)R_{\varphi }^{-k-m}{\left\Vert{\omega
 }\right\Vert}_{L^{r}(B)},$\ \par 
The proof is complete. $\blacksquare $\ \par 

\section{Vitali covering.}

\begin{lem}
Let ${\mathcal{F}}$ be a collection of balls $\lbrace B(x,r(x))\rbrace
 $ in a metric space, with $\forall B(x,r(x))\in {\mathcal{F}},\
 0<r(x)\leq R.$ There exists a disjoint subcollection ${\mathcal{G}}$
 of ${\mathcal{F}}$ with the following properties:\par 
\quad \quad every ball $B$ in ${\mathcal{F}}$ intersects a ball $C$ in ${\mathcal{G}}$
 and $B\subset 5C.$
\end{lem}
This is a well known lemma, see for instance ~\cite{EvGar92},
 section 1.5.1.\ \par 
\quad Fix $\epsilon >0$ and let $\displaystyle \forall x\in M,\ r(x):=R_{\epsilon
 }(x)/120,\ $where $\displaystyle R_{\epsilon }(x)$ is the $(m,\epsilon
 )$-admissible  radius at $\displaystyle x,$ we built a Vitali
 covering with the collection ${\mathcal{F}}:=\lbrace B(x,r(x))\rbrace
 _{x\in M}.$ The previous lemma gives a disjoint subcollection
 ${\mathcal{G}}$ such that every ball $B$ in ${\mathcal{F}}$
 intersects a ball $C$ in ${\mathcal{G}}$ and we have $\displaystyle
 B\subset 5C.$ We set ${\mathcal{G}}':=\lbrace x_{j}\in M::B(x_{j},r(x_{j}))\in
 {\mathcal{G}}\rbrace $ and ${\mathcal{C}}_{\epsilon }:=\lbrace
 B(x,5r(x)),\ x\in {\mathcal{G}}'\rbrace $: we shall call ${\mathcal{C}}(\epsilon
 )$ the $m,\epsilon $ \textbf{admissible covering} of $\displaystyle
 (M,g).$\ \par 
We shall fix $m\geq 0$ and we omit it in order to ease the notation.\ \par 
\ \par 
\quad Then we have the Proposition 7.3 in~\cite{ellipticEq18}:\ \par 

\begin{prop}
~\label{CF9}Let $M$ be a Riemannian manifold, then the overlap
 of the $\epsilon $ admissible covering ${\mathcal{C}}(\epsilon
 )$ is less than $\displaystyle T=\frac{(1+\epsilon )^{n/2}}{(1-\epsilon
 )^{n/2}}(120)^{n},$ i.e.\par 
\quad \quad \quad $\forall x\in M,\ x\in B(y,5r(y))$ where $B(y,r(y))\in {\mathcal{G}}$
 for at most $T$ such balls.\par 
So we have\par 
\quad \quad \quad $\forall f\in L^{1}(M),\ \sum_{j\in {\mathbb{N}}}{\int_{B_{j}}{\left\vert{f(x)}\right\vert
 dv_{g}(x)}}\leq T{\left\Vert{f}\right\Vert}_{L^{1}(M)}.$
\end{prop}

\section{The threshold.~\label{p32}}

\begin{thm}
~\label{HF10}Let $M$ be a Riemannian manifold. Let, for $\displaystyle
 t\geq 0,\ \omega \in L^{2}(M).$ Then we have a solution $u$
 of the heat equation $\displaystyle \partial _{t}u-\Delta u=0,\
 u(x,0)=\omega (x),$ such that $\displaystyle \forall t\geq 0,\
 u(x,t)\in L^{2}(M)$ with the estimate:\par 
\quad \quad \quad $\displaystyle \forall t\geq 0,\ {\left\Vert{u(\cdot ,t)}\right\Vert}_{L^{2}(M)}\leq
 {\left\Vert{\omega }\right\Vert}_{L^{2}(M)}.$
\end{thm}
\quad Proof.\ \par 
It is well known that the Hodge laplacian is essentially positive
 on $p$-forms in $L^{2}(M),$ so $(e^{-t\Delta })_{t\geq 0}$ is
 a contraction semi-group on $\displaystyle L^{2}(M).$ $\blacksquare $\ \par 

\section{Global results.}
\quad We want to globalise Theorem~\ref{m6} by use of our Vitali covering.\ \par 
First set\ \par 
\quad \quad \quad $\displaystyle {\mathcal{D}}(\epsilon ):=\lbrace x\in M::B(x,R_{\epsilon
 }(x))\in {\mathcal{C}}(\epsilon )\rbrace .$\ \par 

\begin{lem}
~\label{m8}Let $f$ be a $p$-form in $M$ and $\tau \in \lbrack
 1,\infty ).$ Set $w(x):=R_{\epsilon }(x)^{\gamma }$ for a $\gamma
 \in {\mathbb{R}}$ and $B(x):=B(x,R_{\epsilon }(x)/10),$ where
 $R_{\epsilon }(x)$ is the $\epsilon $-admissible radius. For
 $l\geq 0$ we have that:\par 
\quad \quad \quad $\displaystyle \forall \tau \geq 1,\ {\left\Vert{\nabla ^{l}f}\right\Vert}_{L^{\tau
 }(M,\ w)}^{\tau }\simeq \sum_{x\in {\mathcal{D}}(\epsilon )}{R_{\epsilon
 }(x)^{\gamma }{\left\Vert{\nabla ^{l}f}\right\Vert}_{L^{\tau
 }(B(x))}^{\tau }}.$
\end{lem}
\quad Proof.\ \par 
Let $x\in {\mathcal{D}}(\epsilon ),$ this implies that $B(x):=B(x,R_{\epsilon
 }(x)/10)\in {\mathcal{C}}(\epsilon ).$ \ \par 
\quad $\bullet $ First we start with $l=0.$ We shall deal with the
 function $\left\vert{f}\right\vert .$\ \par 
\quad For any $y\in B(x)$ we set $R(y)\ :=R_{\epsilon }\ (y).$ We have,
 because ${\mathcal{C}}(\epsilon )$ is a covering of $M$:\ \par 
\quad \quad \quad $\displaystyle {\left\Vert{f}\right\Vert}_{L^{\tau }(M,w)}^{\tau
 }:=\int_{M}{\left\vert{f(x)}\right\vert ^{\tau }w(x)dv(x)}\leq
 \sum_{x\in {\mathcal{D}}(\epsilon )}{\int_{B(x)}{\left\vert{f(y)}\right\vert
 ^{\tau }R(y)^{\gamma }}dv(y)}.$\ \par 
We have, by Lemma~\ref{m7}, $\forall y\in B,\ R(y)\leq 2R(x),$ then\ \par 
\quad \quad \quad $\displaystyle \sum_{x\in {\mathcal{D}}(\epsilon )}{\int_{B(x)}{\left\vert{f(y)}\right\vert
 ^{\tau }R(y)^{\gamma }}dv(y)}\leq $\ \par 
\quad \quad \quad \quad \quad \quad \quad \quad \quad $\displaystyle \leq \sum_{x\in {\mathcal{D}}(\epsilon )}{2^{\gamma
 }R(x)^{\gamma }\int_{B(x)}{\left\vert{f(y)}\right\vert ^{\tau
 }}dv(y)}\leq 2^{\gamma }\sum_{x\in {\mathcal{D}}(\epsilon )}{R(x)^{\gamma
 }{\left\Vert{f}\right\Vert}_{L^{\tau }(B(x))}^{\tau }}.$\ \par 
Hence\ \par 
\quad \quad \quad $\displaystyle {\left\Vert{f}\right\Vert}_{L^{\tau }(M,w)}^{\tau
 }\leq 2^{\gamma }\sum_{x\in {\mathcal{D}}(\epsilon )}{R(x)^{\gamma
 }{\left\Vert{f}\right\Vert}_{L^{\tau }(B)}^{\tau }}.$\ \par 
\ \par 
\quad To get the converse inequality we still use Lemma~\ref{m7}: $\forall
 y\in B,\ R(x)\leq 2R(y)$ so we get:\ \par 
\quad \quad \quad $\displaystyle \sum_{x\in {\mathcal{D}}(\epsilon )}{R(x)^{\gamma
 }\int_{B(x)}{\left\vert{f(y)}\right\vert ^{\tau }}dv(y)}\leq
 2^{\gamma }\sum_{x\in {\mathcal{D}}(\epsilon )}{\int_{B(x)}{R(y)^{\gamma
 }\left\vert{f(y)}\right\vert ^{\tau }}dv(y)}.$\ \par 
Now we use the fact that the overlap of ${\mathcal{C}}(\epsilon
 )$ is bounded by $T,$\ \par 
\quad \quad \quad $\displaystyle \sum_{x\in {\mathcal{D}}(\epsilon )}{\int_{B(x)}{R(y)^{\gamma
 }\left\vert{f(y)}\right\vert ^{\tau }}dv(y)}\leq 2^{\gamma }T\int_{M}{R(y)^{\gamma
 }\left\vert{f(y)}\right\vert ^{\tau }}dv(y)=2^{\gamma }T{\left\Vert{f}\right\Vert}_{L^{\tau
 }(M,w)}^{\tau }.$\ \par 
So\ \par 
\quad \quad \quad $\displaystyle \sum_{x\in {\mathcal{D}}(\epsilon )}{R^{\gamma
 }{\left\Vert{f}\right\Vert}_{L^{\tau }(B)}}^{\tau }\leq 2^{\gamma
 }T{\left\Vert{f}\right\Vert}_{L^{\tau }(M,w)}^{\tau }.$\ \par 
\ \par 
\quad $\bullet $ Now let $l\geq 1.$\ \par 
We apply the case $l=0$ to the covariant derivatives of $f.$\ \par 
\quad \quad \quad $\displaystyle \forall \tau \geq 1,\ {\left\Vert{\nabla ^{l}f}\right\Vert}_{L^{\tau
 }(M,w)}^{\tau }\simeq \sum_{x\in {\mathcal{D}}(\epsilon )}{R(x)^{\gamma
 }{\left\Vert{\nabla ^{l}f}\right\Vert}_{L^{\tau }(B(x))}^{\tau }}.$\ \par 
The proof is complete. $\blacksquare $\ \par 
\ \par 
\quad Let $\omega \in L^{2}(M)\cap L^{r}(M)$ and let $u:=e^{-t\Delta
 }\omega $  be the canonical solution of the heat equation given
 by Theorem~\ref{HF10}, i.e.\ \par 
\quad \quad \quad $\displaystyle \forall t\geq 0,\ {\left\Vert{u(\cdot ,t)}\right\Vert}_{L^{2}(M)}\leq
 {\left\Vert{\omega }\right\Vert}_{L^{2}(M)}.$\ \par 
By Lemma~\ref{m8} we get, replacing $f$ by $u$ and $\tau $ by
 $r,$ with $w(x):=R_{\epsilon }(x)^{\gamma }$ and using the covering
 ${\mathcal{C}}(\epsilon )$:\ \par 
\quad \quad \quad \begin{equation}  \forall r\geq 1,\ {\left\Vert{\nabla ^{l}u(\cdot
 ,t)}\right\Vert}_{L^{r}(M,\ w)}^{r}\simeq \sum_{x\in {\mathcal{D}}(\epsilon
 )}{R_{\epsilon }(x)^{\gamma }{\left\Vert{\nabla ^{l}u(\cdot
 ,t)}\right\Vert}_{L^{r}(B(x))}^{r}}.\label{gC18}\end{equation}\ \par 
But Theorem~\ref{m6} tells us with $B:=B(x,R)\in {\mathcal{A}}_{k}(\epsilon
 )$ for $p$-form with $p\geq 1$ and $B:=B(x,R)\in {\mathcal{A}}_{k-1}(\epsilon
 )$  if $p=0$ and with $m=1$ on functions and $m=2$ on $p$-forms
 with $p\geq 1:$\ \par 
\quad \quad \quad $\displaystyle \forall t\in (\delta ,1),\ {\left\Vert{\nabla
 ^{k}u(\cdot ,t)}\right\Vert}_{L^{r}(B(x,R))}\leq c(n,r)\frac{\delta
 }{t-\delta ^{3/2}}R_{}^{-k-m}{\left\Vert{\omega }\right\Vert}_{L^{r}(B)}.$\
 \par 
so, with $l=k$ in~(\ref{gC18}) we get:\ \par 
\quad \quad \quad $\displaystyle \forall r\geq 1,\ {\left\Vert{\nabla ^{k}u(\cdot
 ,t)}\right\Vert}_{L^{r}(M,\ w)}^{r}\leq \sum_{x\in {\mathcal{D}}(\epsilon
 )}{R_{\epsilon }(x)^{\gamma }{\left\Vert{\nabla ^{k}u(\cdot
 ,t)}\right\Vert}_{L^{r}(B(x))}^{r}}\leq $\ \par 
\quad \quad \quad \quad \quad \quad \quad \quad \quad \quad \quad \quad \quad \quad \quad $\leq c(n,r)^{r}(\frac{\delta }{t-\delta ^{3/2}})^{r}\sum_{x\in
 {\mathcal{D}}(\epsilon )}{R_{\epsilon }(x)^{\gamma }R_{\epsilon
 }(x)^{-kr-rm}{\left\Vert{\omega }\right\Vert}_{L^{r}(B(x))}^{r}}\leq $\ \par 
\quad \quad \quad \quad \quad \quad \quad \quad \quad \quad \quad \quad \quad \quad \quad $\leq c(n,r)^{r}(\frac{\delta }{t-\delta ^{3/2}})^{r}{\left\Vert{\omega
 }\right\Vert}_{L^{r}(M,w')}^{r},$\ \par 
Here we have set $w'(x):=R_{\epsilon }(x)^{\gamma -kr-rm}$ and
 $w(x):=R_{\epsilon }(x)^{\gamma }.$ Hence:\ \par 
\quad \quad \quad $\displaystyle \forall t\in (\delta ,1),\ \forall r\geq 1,\ {\left\Vert{\nabla
 ^{k}u(\cdot ,t)}\right\Vert}_{L^{r}(M,\ w)}\leq c(n,r)(\frac{\delta
 }{t-\delta ^{3/2}}){\left\Vert{\omega }\right\Vert}_{L^{r}(M,w')}.$\ \par 
\ \par 
For $t\geq 1,$ we get by Theorem~\ref{m6} for $p\geq 1,$ hence $m=2$:\ \par 
\quad \quad \quad $\displaystyle \forall r\geq 1,\ {\left\Vert{\nabla ^{k}u(\cdot
 ,t)}\right\Vert}_{L^{r}(B)}\leq c(n,r)R^{-k-2}{\left\Vert{\omega
 }\right\Vert}_{L^{r}(B)},$\ \par 
Exactly the same way as above, we get:\ \par 
\quad \quad \quad $\displaystyle \forall t\geq 1,\ \forall r\geq 1,\ {\left\Vert{\nabla
 ^{k}u(\cdot ,t)}\right\Vert}_{L^{r}(M,\ w)}\leq c(n,r){\left\Vert{\omega
 }\right\Vert}_{L^{r}(M,w')}.$\ \par 
Here we also set $w'(x):=R_{\epsilon }(x)^{\gamma -kr-2r}$ and
 $w(x):=R_{\epsilon }(x)^{\gamma }.$\ \par 
\quad For $p=0,$ i.e. for functions, hence $m=1$\!\!\!\! , with $\displaystyle
 w'(x):=R_{\epsilon }(x)^{\gamma -kr-r}$ and $w(x):=R_{\epsilon
 }(x)^{\gamma },$  we get the better result:\ \par 
\quad \quad \quad $\displaystyle \forall t\geq 1,\ \forall r\geq 1,\ {\left\Vert{\nabla
 ^{k}u(\cdot ,t)}\right\Vert}_{L^{r}(M,\ w)}\leq c(n,r)t^{-1/2}{\left\Vert{\omega
 }\right\Vert}_{L^{r}(M,w')}.$\ \par 
Always with $m=1$ if $p=0$ and $m=2$ if $p\geq 1,$ we have: if
 $k\leq 1$ we have $R_{\epsilon }(x):=R_{m,\epsilon }(x),$ i.e.
 this is the admissible radius for the $(m,\epsilon )$-admissible
 balls. If $k\geq 2,$ then $R_{\epsilon }(x)$ is the admissible
 radius for the $(k,\epsilon )$-admissible balls for $p\geq 1$
 and for the $(k-1,\epsilon )$-admissible balls for $p=0.$ \ \par 
\ \par 
\quad Now we choose, for instance, $\gamma =kr+rm$ and we get, with
 $w(x):=R_{\epsilon }(x)^{kr+rm}$:\ \par 
\quad \quad \quad $\displaystyle \forall t\in (\delta ,1),\ \forall r\geq 1,\ {\left\Vert{\nabla
 ^{k}u(\cdot ,t)}\right\Vert}_{L^{r}(M,\ w)}\leq c(n,r)(\frac{\delta
 }{t-\delta ^{3/2}}){\left\Vert{\omega }\right\Vert}_{L^{r}(M)}.$\ \par 
And, for $p\geq 1$:\ \par 
\quad \quad \quad $\displaystyle \forall t\geq 1,\ \forall r\geq 1,\ {\left\Vert{\nabla
 ^{k}u(\cdot ,t)}\right\Vert}_{L^{r}(M,\ w)}\leq c(n,r){\left\Vert{\omega
 }\right\Vert}_{L^{r}(M)}.$\ \par 
And for functions:\ \par 
\quad \quad \quad $\displaystyle \forall t\geq 1,\ \forall r\geq 1,\ {\left\Vert{\nabla
 ^{k}u(\cdot ,t)}\right\Vert}_{L^{r}(M,\ w)}\leq c(n,r)t^{-1/2}{\left\Vert{\omega
 }\right\Vert}_{L^{r}(M)}.$\ \par 
\quad For $r=\infty ,$ the passage form the local estimates to the
 global ones are obvious, so we proved:\ \par 

\begin{thm}
~\label{CF7}Let $M$ be a Riemannian manifold. Let $r\in \lbrack
 1,\infty \rbrack $ and $\omega \in L^{r}(M)\cap L^{2}(M).$ We
 have, with $u:=e^{t\Delta }\omega $ the canonical solution of
 the heat equation, for any $k\geq 0$:\par 
\quad \quad \quad $\displaystyle \forall t\in (\delta ,1),\ {\left\Vert{\nabla
 ^{k}u(\cdot ,t)}\right\Vert}_{L^{r}(M,\ w)}\leq c(n,r)(\frac{\delta
 }{t-\delta ^{3/2}}){\left\Vert{\omega }\right\Vert}_{L^{r}(M)}.$\par 
Now for $p$-forms with $p\geq 1$:\par 
\quad \quad \quad $\displaystyle \forall t\geq 1,\ {\left\Vert{\nabla ^{k}u(\cdot
 ,t)}\right\Vert}_{L^{r}(M,\ w)}\leq c(n,r){\left\Vert{\omega
 }\right\Vert}_{L^{r}(M)},$\par 
and for functions:\par 
\quad \quad \quad $\displaystyle \forall t\geq 1,\ {\left\Vert{\nabla ^{k}u(\cdot
 ,t)}\right\Vert}_{L^{r}(M,\ w)}\leq c(n,r)t^{-1/2}{\left\Vert{\omega
 }\right\Vert}_{L^{r}(M)}.$\par 
where $w(x):=R_{\epsilon }(x)^{kr+mr}.$\par 
\quad If $k\leq 1$ we have $R_{\epsilon }(x):=R_{m,\epsilon }(x)$ is
 the admissible radius for the $(m,\epsilon )$-admissible balls,
 with $m=1$ if $p=0$ and $m=2$ if $p\geq 1.$  If $k\geq 2,$ then
 $R_{\epsilon }(x)$ is the admissible radius for the $(k,\epsilon
 )$-admissible balls for $p\geq 1$ and for the $(k-1,\epsilon
 )$-admissible balls for $p=0.$
\end{thm}
\quad In particular, making $k=0$ and $k=1$ and with the same conditions
 as above on the admissible balls:\ \par 

\begin{cor}
We have with $r\in \lbrack 1,\infty \rbrack $ and $w(x):=R_{\epsilon
 }(x)^{rm}$\!\!\!\! , with $m=1$ if $p=0$ and $m=2$ if $p\geq 1.$\par 
\quad \quad \quad $\displaystyle \forall t\in (\delta ,1),\ {\left\Vert{u(\cdot
 ,t)}\right\Vert}_{L^{r}(M,w)}\leq c(n,r)(\frac{\delta }{t-\delta
 ^{3/2}}){\left\Vert{\omega }\right\Vert}_{L^{r}(M)}$\par 
and\par 
\quad \quad \quad $\displaystyle \forall t\geq 1,\ {\left\Vert{u(\cdot ,t)}\right\Vert}_{L^{r}(M,w)}\leq
 c(n,r){\left\Vert{\omega }\right\Vert}_{L^{r}(M)}.$\par 
\quad For the gradient estimate, with $w(x):=R_{\epsilon }(x)^{r+rm}$ this time:\par 
\quad \quad \quad $\displaystyle \forall t\in (\delta ,1),\ {\left\Vert{\nabla
 u(\cdot ,t)}\right\Vert}_{L^{r}(M,w)}\leq c(n,r)\frac{\delta
 }{t-\delta ^{3/2}}{\left\Vert{\omega }\right\Vert}_{L^{r}(M)}.$\par 
Now for $p$-forms with $p\geq 1$:\par 
\quad \quad \quad $\displaystyle \forall t\geq 1,\ {\left\Vert{\nabla u(\cdot ,t)}\right\Vert}_{L^{r}(M,w)}\leq
 c(n,r){\left\Vert{\omega }\right\Vert}_{L^{r}(M)}$\par 
and for functions:\par 
\quad \quad \quad $\displaystyle \forall t\geq 1,\ {\left\Vert{\nabla u(\cdot ,t)}\right\Vert}_{L^{r}(M,w)}\leq
 c(n,r)t^{-1/2}{\left\Vert{\omega }\right\Vert}_{L^{r}(M)}.$
\end{cor}

\section{Classical estimates.}
\quad We shall give some examples where we have classical estimates
 using that for any $\displaystyle x\in M,$ we have $\displaystyle
 R_{\epsilon }(x)\geq \eta ,$ via~\cite[Corollary, p. 7] {HebeyHerzlich97}
  (see also Theorem 1.3 in the book by Hebey~\cite{Hebey96}):\ \par 

\begin{cor}
~\label{CF8}Let $M$ be a Riemannian manifold. Let $k\geq 1$;
 if we have the injectivity radius $\displaystyle r_{inj}(x)\geq
 i>0$ and  $\displaystyle \forall j\leq k-1,\ \left\vert{\nabla
 ^{j}Rc_{(M,g)}(x)}\right\vert \leq c$ for all $\displaystyle
 x\in M,$ then there exists a constant $\eta >0,$ depending only
 on $\displaystyle n,\epsilon ,i,k$ and  $\displaystyle c,$ such
 that: $\displaystyle \forall x\in M,\ R_{k,\epsilon }(x)\geq \eta .$\par 
\quad For $k=0,$ if we have the injectivity radius $\displaystyle r_{inj}(x)\geq
 i>0$ and  $\displaystyle Rc_{(M,g)}(x)\geq \lambda g_{x}$ for
 some $\lambda \in {\mathbb{R}}$ and for all $\displaystyle x\in
 M,$ then there exists a constant $\eta >0,$ depending only on
 $\displaystyle n,\epsilon ,i,$ and  $\displaystyle \lambda ,$
 such that: $\displaystyle \forall x\in M,\ R_{0,\epsilon }(x)\geq \eta .$
\end{cor}
\quad Proof.\ \par 
The Theorem of Hebey and Herzlich gives that, under these hypotheses,
 for any $\alpha \in (0,1)$ there exists a constant $\eta >0,$
 depending only on $\displaystyle n,\epsilon ,i,k,\alpha $ and
  $\displaystyle c,$ such that:\ \par 
\quad \quad \quad $\displaystyle \forall x\in M,\ r_{H}(1+\epsilon ,k,\alpha )(x)\geq
 \eta .$\ \par 
So taking our definition with a harmonic coordinates patch, we have that:\ \par 
\quad \quad \quad $\displaystyle R_{k,\epsilon }(x)\geq r_{H}(1+\epsilon ,k,\alpha )(x).$\ \par 
So, a fortiori, this is true when we take the sup for $R_{k,\epsilon
 }(x)$ on \emph{any} smooth coordinates patch, not necessarily
 harmonic coordinates one. $\blacksquare $\ \par 
\ \par 
Then we get our "classical estimates":\ \par 

\begin{thm}
~\label{GC13}Let $M$ be a Riemannian manifold. Let $r\in \lbrack
 1,\infty \rbrack $ and $\omega \in L^{r}(M)\cap L^{2}(M).$ For
 $k=0,\ 1$ suppose that $(M,g)$ has $1$-order weak bounded geometry
 for $p$-forms with $p\geq 1$ and $0$-order weak bounded geometry
 for functions. For $k\geq 2$ suppose that $M$ has $k$-order
 weak bounded geometry for $p$-forms with $p\geq 1$ and $k-1$-order
 weak bounded geometry for functions.\par 
\quad Then the canonical solution $u:=e^{t\Delta }\omega $ of the heat
 equation is such that, for any $k\geq 0$ and with $\eta =\eta
 (n,\epsilon ,i,k)$ given by the Corollary~\ref{CF8} of Hebey and Herzlich:\par 
\quad \quad \quad $\displaystyle \forall t\in (\delta ,1),\ {\left\Vert{\nabla
 ^{k}u(\cdot ,t)}\right\Vert}_{L^{r}(M)}\leq c(n,r,\eta )(\frac{\delta
 }{t-\delta ^{3/2}}){\left\Vert{\omega }\right\Vert}_{L^{r}(M)}.$\par 
And for $p$-forms with $p\geq 1$:\par 
\quad \quad \quad $\displaystyle \forall t\geq 1,\ {\left\Vert{\nabla ^{k}u(\cdot
 ,t)}\right\Vert}_{L^{r}(M)}\leq c(n,r,\eta ){\left\Vert{\omega
 }\right\Vert}_{L^{r}(M)},$\par 
and for functions:\par 
\quad \quad \quad $\displaystyle \forall t\geq 1,\ {\left\Vert{\nabla ^{k}u(\cdot
 ,t)}\right\Vert}_{L^{r}(M)}\leq c(n,r,\eta )t^{-1/2}{\left\Vert{\omega
 }\right\Vert}_{L^{r}(M)}.$
\end{thm}
\quad Proof.\ \par 
We apply Theorem~\ref{CF7} together with Corollary~\ref{CF8}
 to have that there exists $\eta >0$ such that for any $x\in
 M,$  we get $\eta \leq R_{\epsilon }(x)\leq 1.$ Hence:\ \par 
\quad \quad \quad $\displaystyle {\left\Vert{\nabla ^{k}u(\cdot ,t)}\right\Vert}_{L^{r}(M,\
 w)}^{r}:=\int_{M}{\left\vert{\nabla ^{k}u(x,t)}\right\vert ^{r}R_{\epsilon
 }(x)^{r+rm}dx}\geq \eta ^{r(1+m)}{\left\Vert{\nabla ^{k}u(\cdot
 ,t)}\right\Vert}_{L^{r}(M)}^{r}.$\ \par 
So Theorem~\ref{CF7} ends the proof with $m$ as in Theorem~\ref{CF7}
 and with the constant $\displaystyle c(n,r,\eta ):=c(n,r)\eta
 ^{-(1+m)}.$\ \par 
Hence we can forget the weight. $\blacksquare $\ \par 

\section{Appendix 1}

\begin{prop}
Let $r\in \lbrack 1,\infty \rbrack .$ Let $x\in M,$ a Riemannian
 manifold. With $B:=B(x,R)$ a $(m,\epsilon )$-admissible ball
 in $M,$ and any $\delta \in (0,1),$ there is a $\epsilon (\delta
 )>0$ such that for any $\epsilon \leq \epsilon (\delta ),$ if
 $\omega $ is a $p$-form in $\displaystyle L^{r}(B),$ then the
 $p$-form $u_{\varphi }:=e^{-t\Delta _{\varphi }}\omega _{\varphi
 }$ verifies, in ${\mathbb{R}}^{n}$\!\!\!\! , with $\forall \gamma
 \in {\mathbb{N}}^{n},$ $l:=\left\vert{\gamma }\right\vert /2$\par 
\quad \quad \quad $\displaystyle \forall t\in (\delta ,1),\ {\left\Vert{\partial
 ^{\gamma }u_{\varphi }}\right\Vert}_{L^{r}(B_{\varphi })}\leq
 c(n,r)\frac{\delta }{t-\delta ^{1+l}}R_{\varphi }^{-m}{\left\Vert{\omega
 _{\varphi }}\right\Vert}_{L^{r}(B_{\varphi })}.$\par 
And\par 
\quad \quad \quad $\displaystyle \forall t\geq 1,\ {\left\Vert{\partial ^{\gamma
 }u_{\varphi }}\right\Vert}_{L^{r}(B_{\varphi })}\leq c(n,r,\delta
 )t^{-l}R_{\varphi }^{-m}{\left\Vert{\omega _{\varphi }}\right\Vert}_{L^{r}(B_{\varphi
 })}$\par 
with $m=1$ if $p=0$ and $m=2$ if $p\geq 1.$ And also $B_{\varphi
 }=\varphi (B),\ \omega _{\varphi }=\varphi ^{*}\omega $ etc...
\end{prop}
\quad Proof.\ \par 
First we work with the first convolution in~\ref{CF5}.\ \par 
\quad Because we stay in ${\mathbb{R}}^{n}$ and for easing the notation,
 we forget the subscript $\varphi ,$ so we write $u$ for $u_{\varphi
 },\ \omega $ for $\omega _{\varphi }$ etc.\ \par 
\quad Set $\displaystyle Y:=\Delta _{\varphi }-\Delta ,$ by Lemma~\ref{gC16},
 we have:\ \par 
\quad \quad \quad $Yf=\sum_{i,j}{a_{ij}\partial ^{2}_{ij}f}+\sum_{i}{b_{i}\partial _{i}f},$\ \par 
and, for $p=0,$ in $\varphi (B)$ with $B\in {\mathcal{A}}_{1}(\epsilon
 )$:\ \par 
(A)           $\sum_{i,j}{\left\vert{a_{ij}}\right\vert }\leq
 \epsilon ,\ \sum_{i}{\left\vert{b_{i}}\right\vert }\leq C\epsilon
 R^{-1},$\ \par 
and for $p\geq 1$ with $B\in {\mathcal{A}}_{2}(\epsilon )$:\ \par 
(B)           $\sum_{i,j}{\left\vert{a_{ij}}\right\vert }\leq
 \epsilon ,\ \sum_{i}{\left\vert{b_{i}}\right\vert }\leq C\epsilon
 R^{-2}.$\ \par 
\ \par 
\quad Because the (Hodge) laplacian in ${\mathbb{R}}^{n}$ acts on $p$-forms
 componentwise, we fix $t$ and we have:\ \par 
\quad \quad \quad $\displaystyle (Ye^{-t\Delta }\omega )(y,t)=Y\lbrack \int{\Phi
 (y-z,t)\omega (z)dz}\rbrack $\ \par 
with $\Phi $ the heat kernel in ${\mathbb{R}}^{n}.$\ \par 
Set $\displaystyle Y_{kl}:=\frac{\partial ^{2}}{\partial y_{k}\partial
 y_{l}}$ and:\ \par 
\quad \quad \quad $\displaystyle \psi (y,t):=Y_{kl}\int{(\Phi (y-z,t))\omega (z)dz}.$\ \par 
So, again because $\Phi $ is the heat kernel in ${\mathbb{R}}^{n},$\ \par 
\quad \quad \quad $\displaystyle \partial ^{\gamma }_{x}(e^{-t\Delta }\ast (Y_{kl}e^{-t\Delta
 }\omega ))=\partial ^{\gamma }_{x}(e^{-t\Delta }\ast (\partial
 ^{2}_{y_{k}y_{l}}e^{-t\Delta }\omega ))(x,t)=\int{\partial ^{\gamma
 }_{x}\Phi (x-y,t)\psi (y,t)dy}.$\ \par 
Now recall that $\displaystyle \psi (y,t):=\int{\partial ^{2}_{y_{k}y_{l}}(\Phi
 (y-z,t))\omega (z)dz}$ then, extending $\omega $ by $0$ outside
 $\varphi (B),$ by Corollary~\ref{HF1} in the Appendix 2 and
 the inequalities (A) and (B) above:\ \par 
\quad \quad \quad \begin{equation} {\left\Vert{ \psi (\cdot ,t)}\right\Vert}_{L^{r}(B)}\leq
 {\left\Vert{\psi (\cdot ,t)}\right\Vert}_{L^{r}({\mathbb{R}}^{n})}\leq
 \epsilon c(n,r)t^{-1}{\left\Vert{\omega }\right\Vert}_{L^{r}({\mathbb{R}}^{n})}=\epsilon
 c(n,r)t^{-1}{\left\Vert{\omega }\right\Vert}_{L^{r}(B)}.\label{nC30}\end{equation}\
 \par 
With\ \par 
\quad \quad \quad $\displaystyle \theta (x,t):=\partial ^{\gamma }_{x}(e^{-t\Delta
 }\psi )(x,t)=\int{\partial ^{\gamma }_{x}\Phi (x-y,t)\psi (y,t)dy},$ \ \par 
by Proposition~\ref{HH2}, with $\displaystyle u:=\int_{{\mathbb{R}}^{n}}{f(y)\Phi
 (x-y,t)dy}$ and setting $l:=\left\vert{\gamma }\right\vert /2,$\ \par 
\quad \quad \quad $\displaystyle {\left\Vert{\partial ^{\gamma }u(\cdot ,t)}\right\Vert}_{L^{r}({\mathbb{R}}^{n})}\leq
 \epsilon c(n,r)t^{-l}{\left\Vert{f}\right\Vert}_{L^{r}({\mathbb{R}}^{n})},$\
 \par 
so\ \par 
\quad \quad \quad $\displaystyle {\left\Vert{\theta (\cdot ,t)}\right\Vert}_{L^{r}(B)}\leq
 \epsilon c(n,r)t^{-l}{\left\Vert{\psi (\cdot ,t)}\right\Vert}_{L^{r}(B)}.$\
 \par 
But by~(\ref{nC30}):\ \par 
\quad \quad \quad $\displaystyle {\left\Vert{\psi (\cdot ,t)}\right\Vert}_{L^{r}(B)}\leq
 \epsilon c(n,r)t^{-1}{\left\Vert{\omega }\right\Vert}_{L^{r}({\mathbb{R}}^{n})}.$\
 \par 
So\ \par 
\quad \quad \quad $\displaystyle \forall t>0,\ {\left\Vert{\theta (\cdot ,t)}\right\Vert}_{L^{r}(B)}\leq
 \epsilon c(n,r)t^{-1-l}{\left\Vert{\omega }\right\Vert}_{L^{r}(B)}.$\ \par 
\quad For $Y_{k}:=\partial _{k}$ a first order derivative, we get the same way:\ \par 
\quad \quad \quad $\displaystyle \forall t>0,\ {\left\Vert{\partial ^{\gamma }_{x}(e^{-t\Delta
 }\ast (Y_{k}e^{-t\Delta }\omega ))}\right\Vert}_{L^{r}(B)}\leq
 \epsilon c(n,r)t^{-(1+l)}{\left\Vert{\omega }\right\Vert}_{L^{r}(B)}.$\ \par 
\quad Now using the complete form of $Y,$ we get, for $t<1,$ because
 $t^{-1/2}<t^{-1},$\ \par 
\quad \quad \quad $\displaystyle {\left\Vert{\partial ^{\gamma }(e^{-t\Delta }\ast
 (Ye^{-t\Delta }\omega ))(\cdot ,t)}\right\Vert}_{L^{r}(B)}\leq
 {\left\Vert{\sum_{kl}{\partial ^{\gamma }(e^{-t\Delta }\ast
 (a_{kl}\partial ^{2}_{kl}e^{-t\Delta }\omega )))(\cdot ,t)}}\right\Vert}_{L^{r}(B)}+$\
 \par 
\quad \quad \quad \quad \quad \quad \quad \quad \quad \quad \quad \quad \quad \quad \quad \quad \quad \quad \quad $\displaystyle +{\left\Vert{\sum_{k}{\partial ^{\gamma }(e^{-t\Delta
 }\ast (b_{k}\partial _{k}e^{-t\Delta }\omega )))(\cdot ,t)}}\right\Vert}_{L^{r}(B)}\leq
 $\ \par 
\quad \quad \quad \quad \quad \quad \quad \quad \quad \quad \quad \quad \quad \quad \quad \quad \quad \quad \quad \quad $\displaystyle \leq \epsilon (1+R^{-m})c(n,r)t^{-1-l}{\left\Vert{\omega
 }\right\Vert}_{L^{r}(B)}.$\ \par 
\quad And for $t\geq 1,$\ \par 
\quad \quad \quad $\displaystyle {\left\Vert{\partial ^{\gamma }(e^{-t\Delta }\ast
 (Ye^{-t\Delta }\omega ))(\cdot ,t)}\right\Vert}_{L^{r}(B)}\leq
 \epsilon (1+R^{-m})c(n,r)t^{-(1+l)}{\left\Vert{\omega }\right\Vert}_{L^{r}(B)},$\
 \par 
still with $m=1$ in the case of functions, $p=0,$ and $m=2$ in
 the case of $p$-forms, $p\geq 1.$\ \par 
So we have, because $R\leq 1$:\ \par 
\quad \quad \quad $\displaystyle {\left\Vert{\partial ^{\gamma }(e^{-t\Delta }\ast
 (Ye^{-t\Delta }\omega ))(\cdot ,t)}\right\Vert}_{L^{r}(B)}\leq
 \epsilon c(n,r)t^{-\beta }R^{-m}{\left\Vert{\omega }\right\Vert}_{L^{r}(B)}$\
 \par 
with $\beta =1+\left\vert{\gamma }\right\vert /2$ for $0<t<1$
 and $\beta =(1+\left\vert{\gamma }\right\vert )/2$ for $t\geq 1.$\ \par 
\quad We have to treat $\displaystyle e^{-t\Delta }*(Ye^{-t\Delta })^{*j},\
 j\geq 2.$\ \par 
Fix any $\delta \in (0,1)$ and choose $\epsilon \leq \epsilon
 (\delta ,l)$ such that $\displaystyle \epsilon c(n,r)\leq \delta
 ^{1+l}.$\ \par 
\ \par 
\emph{      Now on we shall always suppose that the }$\epsilon
 $\emph{ appearing in our }$(m,\epsilon )$\emph{-admissible ball
 is less than }$\epsilon (\delta ,l).$\ \par 
\ \par 
Then we have $t\in (\delta ,1)\Rightarrow \epsilon c(n,r)t^{-1-l}<1$
 with $2l=\left\vert{\gamma }\right\vert $ and we have for $j=2:$\ \par 
\quad \quad \quad $\displaystyle {\left\Vert{\partial ^{\gamma }(e^{-t\Delta }\ast
 (Ye^{-t\Delta })^{\ast 2})\omega (\cdot ,t))}\right\Vert}_{L^{r}(B)}={\left\Vert{\partial
 ^{\gamma }(e^{-t\Delta }\ast (Ye^{-t\Delta })\ast (Ye^{-t\Delta
 })^{}\omega (\cdot ,t))}\right\Vert}_{L^{r}(B)}.$\ \par 
Set $\mu :=(Ye^{-t\Delta })\omega (\cdot ,t))$ then\ \par 
\quad $\displaystyle {\left\Vert{\partial ^{\gamma }(e^{-t\Delta }\ast
 (Ye^{-t\Delta })^{\ast 2})\omega (\cdot ,t))}\right\Vert}_{L^{r}(B)}={\left\Vert{\partial
 ^{\gamma }(e^{-t\Delta }\ast (Ye^{-t\Delta })\mu (\cdot ,t))}\right\Vert}_{L^{r}(B)}\leq
 \epsilon c(n,r)t^{-\beta }R^{-m}{\left\Vert{\mu }\right\Vert}_{L^{r}(B)}.$\
 \par 
But, by equation~(\ref{nC30}):\ \par 
\quad \quad \quad ${\left\Vert{\mu }\right\Vert}_{L^{r}(B)}={\left\Vert{(Ye^{-t\Delta
 })\omega (\cdot ,t))}\right\Vert}_{L^{r}(B)}\leq \epsilon c(n,r)t^{-1}{\left\Vert{\omega
 }\right\Vert}_{L^{r}(B)}$\ \par 
so\ \par 
\quad \quad \quad $\displaystyle {\left\Vert{\partial ^{\gamma }(e^{-t\Delta }\ast
 (Ye^{-t\Delta })^{\ast 2})\omega (\cdot ,t))}\right\Vert}_{L^{r}(B)}\leq
 \epsilon ^{2}c(n,r)^{2}t^{-1-\beta }R^{-m}{\left\Vert{\omega
 }\right\Vert}_{L^{r}(B)}.$\ \par 
Again\ \par 
\quad \quad \quad $\displaystyle \partial ^{\gamma }(e^{-t\Delta }\ast (Ye^{-t\Delta
 })^{\ast 3})\omega (\cdot ,t))=\partial ^{\gamma }(e^{-t\Delta
 }\ast (Ye^{-t\Delta })^{\ast 2})\mu (\cdot ,t))$\ \par 
and\ \par 
\quad \quad \quad $\displaystyle {\left\Vert{\partial ^{\gamma }(e^{-t\Delta }\ast
 (Ye^{-t\Delta })^{\ast 3})\omega (\cdot ,t))}\right\Vert}_{L^{r}(B)}\leq
 $\ \par 
\quad \quad \quad \quad \quad \quad \quad $\displaystyle \epsilon ^{2}c(n,r)^{2}t^{-1-\beta }R^{-m}{\left\Vert{\mu
 }\right\Vert}_{L^{r}(B)}\leq \epsilon ^{3}c(n,r)^{3}t^{-2-\beta
 }R^{-m}{\left\Vert{\omega }\right\Vert}_{L^{r}(B)}.$\ \par 
\quad So by induction\ \par 
\quad \quad \quad $\displaystyle {\left\Vert{\partial ^{\gamma }(e^{-t\Delta }\ast
 (Ye^{-t\Delta })^{\ast j})\omega (\cdot ,t))}\right\Vert}_{L^{r}(B)}\leq
 \epsilon ^{j}c(n,r)^{j}t^{-j+1-\beta }R^{-m}{\left\Vert{\omega
 }\right\Vert}_{L^{r}(B)}.$\ \par 
We have, for $\epsilon c(n,r)t^{-1}<1,$\ \par 
\quad \quad \quad $\displaystyle \sum_{j=1}^{\infty }{\epsilon ^{j}c(n,r)^{j}t^{-j}}=\frac{2\epsilon
 c(n,r)}{t-\epsilon c(n,r)}.\ \ \ \ \ (C)$\ \par 
Hence the series\ \par 
\quad \quad \quad $\displaystyle {\left\Vert{\partial ^{\gamma }(e^{-t\Delta }\ast
 \sum_{j=1}^{\infty }{(-1)^{j}(Ye^{-t\Delta })^{\ast j}})\omega
 (\cdot ,t)}\right\Vert}_{L^{r}(B)}\leq (\sum_{j=1}^{\infty }{\epsilon
 ^{j}c(n,r)^{j}t^{-j+1-\beta }})R^{-m}{\left\Vert{\omega }\right\Vert}_{L^{r}(B)}$\
 \par 
converges for $\delta <t<1$ because $\displaystyle \epsilon c(n,r)\leq
 \delta ^{1+l}<1$ and we get:\ \par 
\quad \quad \quad $\displaystyle \sum_{j=1}^{\infty }{\epsilon ^{j}c(n,r)^{j}t^{-j}}=\frac{2\epsilon
 c(n,r)}{t-\epsilon c(n,r)}\leq \frac{2\delta ^{1+l}}{t-\delta ^{1+l}}.$\ \par 
Hence, because $\beta =1+l$ if $t<1,$ we get, with $t^{1-\beta }=t^{-l}$:\ \par 
\quad \quad \quad $\displaystyle \forall t\in (\delta ,1),\ {\left\Vert{\partial
 ^{\gamma }(e^{-t\Delta }\ast \sum_{j=1}^{\infty }{(-1)^{j}(Ye^{-t\Delta
 })^{\ast j}})\omega (\cdot ,t)}\right\Vert}_{L^{r}(B)}\leq 2t^{-l}R^{-m}\frac{\delta
 ^{1+l}}{t-\delta ^{1+l}}{\left\Vert{\omega }\right\Vert}_{L^{r}(B)}.$\ \par 
\ \par 
If $t\geq 1$ we have using $(C)$\ \par 
\quad \quad \quad $\displaystyle (\sum_{j=1}^{\infty }{\epsilon ^{j}c(n,r)^{j}t^{-j+1-\beta
 }})R^{-m}{\left\Vert{\omega }\right\Vert}_{L^{r}(B)}\leq 2t^{-l}R^{-m}\frac{\epsilon
 c(n,r)}{t-\epsilon c(n,r)}\leq 2t^{-l}R^{-m}\frac{\delta ^{1+l}}{t-\delta
 ^{1+l}}{\left\Vert{\omega }\right\Vert}_{L^{r}(B)}.$\ \par 
Hence again\ \par 
\quad \quad \quad $\displaystyle \forall t\geq 1,\ {\left\Vert{\partial ^{\gamma
 }(e^{-t\Delta }\ast \sum_{j=1}^{\infty }{(-1)^{j}(Ye^{-t\Delta
 })^{\ast j}})\omega (\cdot ,t)}\right\Vert}_{L^{r}(B)}\leq 2t^{-l}R^{-m}\frac{\delta
 ^{1+l}}{t-\delta ^{1+l}}{\left\Vert{\omega }\right\Vert}_{L^{r}(B)}.$\ \par 
The formula~(\ref{CF5}) gives:\ \par 
\quad \quad \quad $\displaystyle e^{-t\Delta _{\varphi }}=e^{-t\Delta }+e^{-t\Delta
 }\ast \sum_{j=1}^{\infty }{(-1)^{j}(Ye^{-t\Delta })^{\ast j}},$\ \par 
it remains to add the first term in the right hand side:\ \par 
\quad \quad \quad $\displaystyle {\left\Vert{\partial ^{\gamma }e^{-t\Delta _{\varphi
 }}\omega }\right\Vert}_{L^{r}(B)}\leq {\left\Vert{\partial ^{\gamma
 }e^{-t\Delta }\omega }\right\Vert}_{L^{r}(B)}+\ {\left\Vert{\partial
 ^{\gamma }(e^{-t\Delta }\ast \sum_{j=1}^{\infty }{(-1)^{j}(Ye^{-t\Delta
 })^{\ast j}})\omega }\right\Vert}_{L^{r}(B)}.$\ \par 
Using Corollary~\ref{HF1}:\ \par 
\quad \quad \quad $\displaystyle {\left\Vert{\partial ^{\gamma }e^{-t\Delta }\omega
 }\right\Vert}_{L^{r}({\mathbb{R}}^{n})}\leq c(n,r)t^{-l}{\left\Vert{\omega
 }\right\Vert}_{L^{r}({\mathbb{R}}^{n})},$\ \par 
so adding we get\ \par 
\quad \quad \quad $\displaystyle \forall t\in (\delta ,1),\ {\left\Vert{\partial
 ^{\gamma }e^{-t\Delta _{\varphi }}\omega }\right\Vert}_{L^{r}(B)}\leq
 c(n,r)t^{-l}R^{-m}{\left\Vert{\omega }\right\Vert}_{L^{r}(B)}+2t^{-l}R^{-m}\frac{\delta
 ^{1+l}}{t-\delta ^{1+l}}{\left\Vert{\omega }\right\Vert}_{L^{r}(B)},$\ \par 
hence\ \par 
\quad \quad \quad $\displaystyle \forall t\in (\delta ,1),\ {\left\Vert{\partial
 ^{\gamma }e^{-t\Delta _{\varphi }}\omega }\right\Vert}_{L^{r}(B)}\leq
 t^{-l}R^{-m}{\left\Vert{\omega }\right\Vert}_{L^{r}(B)}(c(n,r)+\frac{2\delta
 ^{1+l}}{t-\delta ^{1+l}})$\ \par 
so with another constant $c(n,r)$\ \par 
\quad \quad \quad $\displaystyle \forall t\in (\delta ,1),\ {\left\Vert{\partial
 ^{\gamma }e^{-t\Delta _{\varphi }}\omega }\right\Vert}_{L^{r}(B)}\leq
 c(n,r)\delta ^{-l}\frac{\delta ^{1+l}}{t-\delta ^{1+l}}R^{-m}{\left\Vert{\omega
 }\right\Vert}_{L^{r}(B)}=$\ \par 
\quad \quad \quad \quad \quad \quad \quad \quad \quad \quad \quad \quad \quad \quad \quad \quad \quad \quad \quad \quad \quad $\displaystyle =c(n,r)\frac{\delta }{t-\delta ^{1+l}}R^{-m}{\left\Vert{\omega
 }\right\Vert}_{L^{r}(B)}.$\ \par 
And, the same way:\ \par 
\quad \quad \quad $\displaystyle \forall t\geq 1,\ {\left\Vert{\partial ^{\gamma
 }e^{-t\Delta _{\varphi }}\omega }\right\Vert}_{L^{r}(B)}\leq
 t^{-l}R^{-m}{\left\Vert{\omega }\right\Vert}_{L^{r}(B)}(c(n,r)+\frac{\delta
 ^{1+l}}{t-\delta ^{1+l}})$\ \par 
so again, still with $l:=\left\vert{\gamma }\right\vert /2$ and
 with another constant $c(n,r)$\ \par 
\quad \quad \quad $\displaystyle \forall t\geq 1,\ {\left\Vert{\partial ^{\gamma
 }e^{-t\Delta _{\varphi }}\omega }\right\Vert}_{L^{r}(B)}\leq
 c(n,r)\frac{\delta ^{1+l}}{1-\delta ^{1+l}}t^{-l}R^{-m}{\left\Vert{\omega
 }\right\Vert}_{L^{r}(B)}.$\ \par 
\quad Recall that we dropped the index $\varphi $ on $u,\ R,\ \omega $ etc...\ \par 
The proof of the Proposition is complete. $\blacksquare $\ \par 

\section{Appendix 2. The heat kernel in ${\mathbb{R}}^{n}.$}
\quad All results here are very well known and they are here essentially
 to fix the notation. See for instance~\cite{Evans98}.\ \par 
\quad We have the heat operator $\displaystyle Du:=\partial _{t}u-\Delta
 u$ and the heat kernel in ${\mathbb{R}}^{n}$:\ \par 
\quad \quad \quad $\displaystyle \Phi (x,t):={\left\lbrace{
\begin{matrix}
{\frac{1}{(4\pi t)^{n/2}}e^{-\frac{\left\vert{x}\right\vert ^{2}}{4t}}}&{\
 x\in {\mathbb{R}}^{n},\ t>0}\cr 
{0}&{\ x\in {\mathbb{R}}^{n},\ t\leq 0}\cr 
\end{matrix}
}\right.}$\ \par 
and an easy computation gives:\ \par 
\quad \quad \quad $\displaystyle \partial _{j}\Phi (x,t)=-\frac{x_{j}}{2t(4\pi
 t)^{n/2}}e^{-\frac{\left\vert{x}\right\vert ^{2}}{4t}}$\ \par 
\quad \quad \quad $\displaystyle \ \partial ^{2}_{j}\Phi (x,t)=(-\frac{1}{2t}+\frac{x_{j}^{2}}{4t^{2}})\frac{1}{(4\pi
 t)^{n/2}}e^{-\frac{\left\vert{x}\right\vert ^{2}}{4t}}$\ \par 
\quad \quad \quad $\displaystyle \partial ^{2}_{jk}\Phi (x,t)=\frac{x_{j}x_{k}}{4t^{2}(4\pi
 t)^{n/2}}e^{-\frac{\left\vert{x}\right\vert ^{2}}{4t}}$\ \par 
\quad \quad \quad \begin{equation}  \ {\left\Vert{\Phi (\cdot ,t)}\right\Vert}_{L^{r}({\mathbb{R}}^{n})}=c_{0}(n,r)\frac{1}{(t)^{\frac{n}{2}(1-\frac{1}{r})}}.\label{HHC6}\end{equation}\
 \par 
and\ \par 
\quad \quad \quad \begin{equation}  \ {\left\Vert{\nabla \Phi (\cdot ,t)}\right\Vert}_{L^{r}({\mathbb{R}}^{n})}\leq
 c_{1}(n,r)\frac{1}{(t)^{\frac{1}{2}+\frac{n}{2}(1-\frac{1}{r})}}.\label{HH0}\end{equation}\
 \par 
and more generally:\ \par 
\quad \quad \quad \begin{equation}  \forall k\in {\mathbb{N}},\ {\left\Vert{\nabla
 ^{k}\Phi (\cdot ,t)}\right\Vert}_{L^{r}({\mathbb{R}}^{n})}\leq
 c_{k}(n,r)\frac{1}{(t)^{\frac{k}{2}+\frac{n}{2}(1-\frac{1}{r})}}.\label{gC14}\end{equation}\
 \par 
\quad These inequalities can be written for $k\in {\mathbb{N}}$: \ \par 
\quad \quad \quad $\displaystyle {\left\Vert{\Phi (\cdot ,t)}\right\Vert}_{W^{k,r}({\mathbb{R}}^{n})}=c_{k}(n,r)\frac{1}{(t)^{\frac{k}{2}+\frac{n}{2}(1-\frac{1}{r})}}.$\
 \par 
\ \par 
\quad We get the following lemma:\ \par 

\begin{lem}
~\label{HH1}We have for $1\leq r\leq s\leq \infty $:\par 
\quad \quad \quad \quad \quad \quad \quad $\displaystyle {\left\Vert{\int_{{\mathbb{R}}^{n}}{f(y)\Phi (x-y,t)dy}}\right\Vert}_{L^{s}({\mathbb{R}}^{n})}\leq
 c(n,r,s)\frac{1}{t^{\frac{n}{2}(\frac{1}{r}-\frac{1}{s})}}{\left\Vert{f}\right\Vert}_{L^{r}({\mathbb{R}}^{n})},$\par
 
where $\displaystyle c(n,r,s)$ is a constant depending only on
 $\displaystyle n,r$ and $\displaystyle s.$\par 
And\par 
\quad \quad \quad \quad \quad \quad \quad $\displaystyle {\left\Vert{\int_{{\mathbb{R}}^{n}}{f(y)\nabla
 \Phi (x-y,t)dy}}\right\Vert}_{L^{s}({\mathbb{R}}^{n})}\leq c(n,r,s)\frac{1}{(t)^{\frac{1}{2}+\frac{n}{2}(\frac{1}{r}-\frac{1}{s})}}{\left\Vert{f}\right\Vert}_{L^{r}({\mathbb{R}}^{n})}.$\par
 
And more generally:\par 
\quad \quad \quad \quad \quad \quad \quad $\displaystyle {\left\Vert{\int_{{\mathbb{R}}^{n}}{f(y)\frac{\partial
 ^{\alpha }}{\partial x_{1}^{\alpha _{1}}...\partial x_{n}^{\alpha
 _{n}}}\Phi (x-y,t)dy}}\right\Vert}_{L^{s}({\mathbb{R}}^{n})}\leq
 c(n,r,s)\frac{1}{(t)^{\frac{\left\vert{\alpha }\right\vert }{2}+\frac{n}{2}(\frac{1}{r}-\frac{1}{s})}}{\left\Vert{f}\right\Vert}_{L^{r}({\mathbb{R}}^{n})},$
\end{lem}

      Proof.\ \par 
The convolution gives:\ \par 
\quad \quad \quad \quad \quad \quad \quad $\displaystyle {\left\Vert{\int_{{\mathbb{R}}^{n}}{f(y)\Phi (x-y,t)dy}}\right\Vert}_{L^{s}({\mathbb{R}}^{n})}\leq
 {\left\Vert{f}\right\Vert}_{L^{r}(M)}{\left\Vert{\Phi (\cdot
 ,t)}\right\Vert}_{L^{u}({\mathbb{R}}^{n})}$\ \par 
with $\displaystyle \ \frac{1}{s}=\frac{1}{r}+\frac{1}{u}-1.$
 Using~(\ref{HHC6}), we get\ \par 
\quad \quad \quad \quad \quad \quad \quad $\displaystyle {\left\Vert{\int_{{\mathbb{R}}^{n}}{f(y)\Phi (x-y,t)dy}}\right\Vert}_{L^{s}({\mathbb{R}}^{n})}\leq
 c_{0}(n,r,u)\frac{1}{(t)^{\frac{n}{2}(1-\frac{1}{u})}}{\left\Vert{f}\right\Vert}_{L^{r}({\mathbb{R}}^{n})}$\
 \par 
hence,\ \par 
\quad \quad \quad \quad \quad \quad \quad $\displaystyle {\left\Vert{\int_{{\mathbb{R}}^{n}}{f(y)\Phi (x-y,t)dy}}\right\Vert}_{L^{s}({\mathbb{R}}^{n})}\leq
 c(n,r,s)\frac{1}{(t)^{\frac{n}{2}(\frac{1}{r}-\frac{1}{s})}}{\left\Vert{f)}\right\Vert}_{L^{r}({\mathbb{R}}^{n})}.$\
 \par 
For the second part we proceed the same way with~(\ref{HH0})
 in place of ~(\ref{HHC6}), to get:\ \par 
\quad \quad \quad \quad \quad \quad \quad $\displaystyle {\left\Vert{\int_{{\mathbb{R}}^{n}}{f(y)\nabla
 \Phi (x-y,t)dy}}\right\Vert}_{L^{s}({\mathbb{R}}^{n})}\leq c(n,r,s)\frac{1}{(t)^{\frac{1}{2}+\frac{n}{2}(1-\frac{1}{s})}}{\left\Vert{f}\right\Vert}_{L^{r}({\mathbb{R}}^{n})}.$\
 \par 
The third part is the same, with~(\ref{gC14}) instead of~(\ref{HH0}),
 which proves the lemma. $\blacksquare $\ \par 

\begin{prop}
~\label{HH2} Let $u(x,t):=\int_{{\mathbb{R}}^{n}}{f(y)\Phi (x-y,t)dy}$
 then, for $1\leq r\leq s\leq \infty ,$ we have, with $k=(k_{1},...,k_{n})$
 and $\left\vert{k}\right\vert =k_{1}+...+k_{n}$:\par 
\quad \quad \quad \quad \quad (i)     $\displaystyle \ {\left\Vert{u(\cdot ,t)}\right\Vert}_{L^{s}({\mathbb{R}}^{n})}\leq
 c(n,r,s)\frac{1}{t^{\frac{n}{2}(\frac{1}{r}-\frac{1}{s})}}{\left\Vert{f}\right\Vert}_{L^{r}({\mathbb{R}}^{n})}.$\par
 
And\par 
\quad \quad \quad \quad \quad (ii)      $\displaystyle \ {\left\Vert{\nabla u(\cdot ,t)}\right\Vert}_{L^{s}({\mathbb{R}}^{n})}\leq
 c(n,r,s)\frac{1}{(t)^{\frac{1}{2}+\frac{n}{2}(\frac{1}{r}-\frac{1}{s})}}{\left\Vert{f}\right\Vert}_{L^{r}({\mathbb{R}}^{n})}.$\par
 
And more generally:\par 
\quad \quad \quad \quad \quad (iii)      $\displaystyle \ {\left\Vert{\frac{\partial ^{\left\vert{k}\right\vert
 }}{\partial x_{1}^{k_{1}}...\partial x_{n}^{k_{n}}}u(\cdot ,t)}\right\Vert}_{L^{s}({\mathbb{R}}^{n})}\leq
 c(n,r,s)\frac{1}{(t)^{\frac{\left\vert{k}\right\vert }{2}+\frac{n}{2}(\frac{1}{r}-\frac{1}{s})}}{\left\Vert{f}\right\Vert}_{L^{r}({\mathbb{R}}^{n})}.$\par
 
\quad This can also be written:\par 
\quad \quad \quad \quad \quad $\displaystyle {\left\Vert{\partial ^{k}(e^{t\Delta }))}\right\Vert}_{L^{r}({\mathbb{R}}^{n})-L^{s}({\mathbb{R}}^{n})}\leq
 c(n,r,s)\frac{1}{(t)^{\frac{\left\vert{k}\right\vert }{2}+\frac{n}{2}(\frac{1}{r}-\frac{1}{s})}}.$
\end{prop}

      Proof.\ \par 
From\ \par 
\quad \quad \quad \quad \quad \quad \quad $\displaystyle u(x,t):=\int_{{\mathbb{R}}^{n}}{f(y)\Phi (x-y,t)dy}$\ \par 
we get\ \par 
\quad \quad \quad \quad \quad \quad \quad $\displaystyle {\left\Vert{u(\cdot ,t)}\right\Vert}_{L^{s}({\mathbb{R}}^{n})}\leq
 {\left\Vert{\int_{{\mathbb{R}}^{n}}{f(y)\Phi (\cdot -y,t)dy}}\right\Vert}_{L^{s}({\mathbb{R}}^{n})},$\
 \par 
hence, by lemma~\ref{HH1},\ \par 
\quad \quad \quad \quad \quad \quad \quad $\displaystyle {\left\Vert{u(\cdot ,t)}\right\Vert}_{L^{s}({\mathbb{R}}^{n})}\leq
 c(n,r,s)\frac{1}{t^{\frac{n}{2}(\frac{1}{r}-\frac{1}{s})}}{\left\Vert{f}\right\Vert}_{L^{r}({\mathbb{R}}^{n})},$\
 \par 
and the (i).\ \par 
\ \par 
\quad The same we get,\ \par 
\quad \quad \quad \quad \quad \quad \quad $\displaystyle \nabla u(x,t):=\int_{{\mathbb{R}}^{n}}{f(y)\nabla
 _{x}\Phi (x-y,t)dy}$\ \par 
hence\ \par 
\quad \quad \quad \quad \quad \quad \quad $\displaystyle {\left\Vert{\nabla u(\cdot ,t)}\right\Vert}_{L^{s}({\mathbb{R}}^{n})}\leq
 {\left\Vert{\int_{{\mathbb{R}}^{n}}{f(y)\nabla \Phi (x-y,t)dy}}\right\Vert}_{L^{s}({\mathbb{R}}^{n})}$\
 \par 
and,  by lemma~\ref{HH1}, we get the (ii).\ \par 
\quad \quad \quad \quad \quad \quad \quad $\displaystyle {\left\Vert{\nabla u(\cdot ,t)}\right\Vert}_{L^{s}({\mathbb{R}}^{n})}\leq
 c(n,r,s)\frac{1}{(t)^{\frac{1}{2}+\frac{n}{2}(\frac{1}{r}-\frac{1}{s})}}{\left\Vert{f}\right\Vert}_{L^{r}({\mathbb{R}}^{n})}.$\
 \par 
\ \par 
\quad The same we get, by Lemma~\ref{HH1},\ \par 
\quad \quad \quad \quad \quad \quad \quad $\displaystyle {\left\Vert{\partial ^{k}u(\cdot ,t)}\right\Vert}_{L^{s}({\mathbb{R}}^{n})}\leq
 c(n,r,s)\frac{1}{(t)^{\frac{\left\vert{k}\right\vert }{2}+\frac{n}{2}(\frac{1}{r}-\frac{1}{s})}}{\left\Vert{f}\right\Vert}_{L^{r}({\mathbb{R}}^{n})}.$\
 \par 
\quad This ends the proof of the proposition. $\blacksquare $\ \par 

\begin{cor}
~\label{HF1}Let $r\in \lbrack 1,\infty \rbrack .$ We have, if
 $u$ is the solution of the heat equation: $\partial _{t}u-\Delta
 u=0,\ u(x,0)=f(x),$ given by $u(^{}x,t):=e^{t\Delta }f=\int_{{\mathbb{R}}^{n}}{f(y)\Phi
 (x-y,t)dy}$\!\!\!\! , with $\displaystyle \partial ^{k}f:=\frac{\partial
 ^{\left\vert{k}\right\vert }}{\partial x_{1}^{k_{1}}...\partial
 x_{n}^{k_{n}}}f$:\par 
\quad \quad \quad $\displaystyle \forall k=(k_{1},...,k_{n})\in {\mathbb{N}}^{n},\
 {\left\Vert{\partial ^{k}u(\cdot ,t)}\right\Vert}_{L^{r}({\mathbb{R}}^{n})}\leq
 c(n,r)\frac{1}{t^{\left\vert{k}\right\vert /2}}{\left\Vert{f}\right\Vert}_{L^{r}({\mathbb{R}}^{n})},$\par
 
or, equivalently:\par 
\quad \quad \quad $\displaystyle \forall k=(k_{1},...,k_{n})\in {\mathbb{N}}^{n},\
 {\left\Vert{\partial ^{k}e^{t\Delta }}\right\Vert}_{L^{r}({\mathbb{R}}^{n})-L^{r}({\mathbb{R}}^{n})}\leq
 c(n,r)\frac{1}{t^{\left\vert{k}\right\vert /2}},$
\end{cor}
\quad Proof.\ \par 
We apply Proposition~\ref{HH2} with $r=s.$ $\blacksquare $\ \par 
\ \par 

\bibliographystyle{/usr/local/texlive/2017/texmf-dist/bibtex/bst/base/plain}

\end{document}